\documentclass[11pt,leqno]{amsart}
\usepackage{amsmath,amsfonts,amssymb,amscd,amsthm,amsbsy,upref}
\numberwithin{equation}{section}
\newtheorem{thm}{Theorem}[section]
\newtheorem{lem}[thm]{Lemma}
\newtheorem{cor}[thm]{Corollary}

\newtheorem{prop}[thm]{Proposition}
\newtheorem*{thmA}{Theorem A}
\newtheorem*{thmB}{Theorem B}
\newtheorem*{corC}{Corollary C}
\newtheorem*{corD}{Corollary D}
\theoremstyle{definition}

\newtheorem{prob}[thm]{Problem}

\newtheorem*{remark}{Remark}
\def\E{{\mathbb E}}
\def\Pr{{\mathbb P}}
\def\nat{{\mathbb N}}
\def\real{{\mathbb R}}
\def\F{{\mathcal F}}
\def\ep{\varepsilon}
\def\L{{\mathcal L}}
\def\Ksim{\mathop{ \ \buildrel K\over\sim\  }\nolimits}
\def\Kpsim{\mathop{ \ \buildrel {K_p}\over\sim\  }\nolimits}
\def\Cpsim{\mathop{ \ \buildrel {C_p}\over\sim\  }\nolimits}
\def\Csim{\mathop{ \ \buildrel {C}\over\sim\  }\nolimits}
\def\Csima#1{\mathop{ \ \buildrel {#1C}\over\sim\  }\nolimits}
\def\Dpsim{\mathop{ \ \buildrel {D_p}\over\sim\  }\nolimits}
\def\wmto{\mathop{ \ \buildrel {\text{wm}}\over\longrightarrow\ }\nolimits}
\def\wto{\mathop{ \ \buildrel {\text{w}}\over\longrightarrow\ }\nolimits}
\def\Qarrow{\mathop{ \ \buildrel {\text{Q}}\over\longrightarrow\ }\nolimits}
\def\cA{\mathcal A}
\def\xb{\overline{x}}
\def\supp{\text{\rm supp}}
\def\epb{\overline{\ep}}
\def\deltab{\overline{\delta}}

\def\trivert{|\!|\!|}
\def\Btrivert{\Big|\!\Big|\!\Big|}

\begin{document}
\title{Small Subspaces of $L_p$}
\author{R. Haydon, E. Odell \and Th. Schlumprecht}
\address{Brasenose College\\ Oxford  OX1 4AJ, U.K.}
\email{richard.haydon@brasenose.oxford.ac.uk}
\address{Department of Mathematics\\ The University of Texas at Austin, Austin, TX 78712-0257}
\email{odell@math.utexas.edu}
\address{Department of Mathematics\\ Texas A\&M University\\
College Station, TX 77843-3368} \email{schlump@math.tamu.edu}
\thanks{Research of the last two authors was partially supported by the National Science Foundation.}
\begin{abstract}
We prove that if $X$ is a subspace of $L_p$ $(2<p<\infty)$, then
either $X$ embeds isomorphically into $\ell_p \oplus \ell_2$ or $X$
contains a subspace $Y,$ which is isomorphic to $\ell_p(\ell_2)$. We also give an intrinsic characterization of
when $X$ embeds into $\ell_p \oplus \ell_2$ in terms of weakly null
trees in $X$ or, equivalently, in terms of the ``infinite asymptotic
game'' played in $X$. This solves problems concerning small subspaces of $L_p$
originating in the 1970's. The techniques used were developed over several  decades,
 the most recent being that of weakly null trees developed in the 2000's.
  
\end{abstract}
\maketitle

\section{Introduction}\label{S:1}

The study of ``small subspaces'' of $L_p$ ($2<p<\infty$) was
initiated by Kadets and Pe{\l}czy\'nski \cite{KP} who proved that if
$X$ is an infinite dimensional  subspace of $L_p$, then either $X$ is
isomorphic to $\ell_2$ and the $L_2$-norm is equivalent to the
$L_p$-norm
 on $X$, or for all $\ep >0$ $X$ contains a
subspace $Y$ which is $1+\ep$-isomorphic to $\ell_p$. In \cite{JO1}
it was shown that if $X$ does not contain an isomorph of $\ell_2$
then $X$ embeds isomorphically into $\ell_p$ (\cite{KW} showed that,
moreover, for all $\ep  >0$, $X$\ \  $1+\ep$-embeds into $\ell_p$).
W.B.~Johnson \cite{J} solved the analogous problem for $X\subseteq
L_p$ (for all $1<p<2$) by proving that $X$ embeds into $\ell_p$ if for some
$K<\infty$ every weakly null sequence in $S_X$, the unit sphere of
$X$, admits a subsequence $K$-equivalent to the unit vector basis of
$\ell_p$.

Using the machinery of \cite{OS1} (see also \cite{OS2})
 and the special nature of $L_p$, these
results were unified in \cite{AO} as: $X\subseteq L_p$
$(1<p<\infty)$ embeds into $\ell_p$ if (and only if) every weakly
null tree in $S_X$ admits a branch equivalent to the unit vector
basis of $\ell_p$.

After $\ell_p$ and $\ell_2$ the next smallest natural subspace of
$L_p$ $(2<p<\infty)$ is $\ell_p\oplus \ell_2$. Indeed if $X\subseteq
L_p$ does not embed into either $\ell_p$ or $\ell_2$, it contains an
isomorph of $\ell_p\oplus \ell_2$. The next small natural subspace
after $\ell_p \oplus \ell_2$ is $\ell_p$ $(\ell_2)$ or, as it is
sometimes denoted, $(\sum \ell_2)_p$. In \cite{JO2} it was shown
that if $X\subseteq L_p$ $(2<p<\infty)$ and $X$ is a quotient of a
subspace of $\ell_p\oplus\ell_2$ then $X$ embeds into
$\ell_p\oplus\ell_2$.

The motivating problem for this paper (and our main result) dates
back to the 1970's. We prove that if $X\subseteq L_p$ $(2<p<\infty)$
and $X$ does not embed into $\ell_p\oplus \ell_2$ then $X$ contains
an isomorph of $\ell_p$ $(\ell_2)$. To solve this we first give an
intrinsic characterization of when $X$ embeds into $\ell_p\oplus
\ell_2$. The terminology is explained in Section~3. We assume that
our space
 $L_p$ is defined over an atomless and separable probability space $(\Omega,\Sigma,\mathbb P)$.
We write $A\Ksim B$ if $A\le KB$ and $B\le KA$. $X$ will always
denote an infinite dimensional Banach space.

\begin{thmA}
Let $X$ be a subspace of $L_p$ $(2<p<\infty)$. Then the following
are equivalent.
\begin{itemize}
\item[a)] $X$ embeds into $\ell_p\oplus\ell_2$.
\item[b)] Every weakly null tree in $S_X$ admits a branch $(x_i)$ satisfying
for some $K$ and all scalars $(a_i)$,
\begin{equation}\label{eqA1}
\Big\|\sum a_i x_i\Big\| \Ksim   \left(\sum |a_i|^p\right)^{1/p}
\vee \Big\|\sum a_i x_i\ \Big\|_2 .
\end{equation}
($\|\cdot\|_2$ denotes the $L_2$-norm)
\item[c)] Every weakly null tree in $S_X$ admits a branch $(x_i)$
satisfying, for some $K$, $(w_i) \subseteq [0,1]$, and all scalars
$(a_i)$
\begin{equation}\label{eqA2}
\Big\|\sum a_i x_i\Big\|  \Ksim  \left(\sum |a_i|^p\right)^{1/p}
\vee \left(\sum |a_i|^2 w_i^2\right)^{1/2}\ .
\end{equation}
\end{itemize}
\end{thmA}

Under any of these conditions the embedding of $X$ into
$\ell_p\oplus\ell_2$ is given by:  producing a blocking $(H_n)$ of
the Haar basis for $L_p$ and $1\le K<\infty$, so that, if $X\ni x =
\sum x_n$, $x_n\in H_n$, then
\begin{equation*}
\|x\| \Ksim \left( \sum \|x_n\|_p^p\right)^{1/p} \vee \left(\sum
\|x_n\|_2^2\right)^{1/2} = \left(\sum \|x_n\|_p^p\right)^{1/p} \vee
\|x\|_2\ .
\end{equation*}
Since $(\sum H_n)_p$ is isomorphic to $ \ell_p$ this suffices.

The next task is to show that if $X$ violates these conditions then
$X$ contains  a complemented subspace isomorphic to $\ell_p$ $(\ell_2)$. We will present two proofs of
this. The first proof  will roughly show that $X$ must contain ``skinny'' uniform
copies of $\ell_2$ and hence contain uniform $\ell_2$'s,
$(X_n)_{n\in\nat}$ for which if $x_n\in S_{X_n}$ then the $x_n$'s
are almost disjointly supported and hence behave like the unit
vector basis of $\ell_p$.  Then an argument due to
 Schechtman   will prove that a subspace of $X$ which is isomorphic
  to $\ell_p(\ell_2)$ contains an isomorphic copy of $\ell_p(\ell_2)$ 
   which is complemented in $L_p$.
  The second proof will lead to    more precise result  using the
random measure machinery of D.~Aldous \cite{A} and the stability
theory of $L_p$ \cite{KM}. For easier reading we will, however, recall
 all relevant definitions and results concerning  random measures and
stability theory.
 We will show that the complemented copy  of $\ell_p(\ell_2)$ 
 is witnessed by {\em stabilized $\ell_2$ sequences } living on
  almost disjoint supports, meaning that the joint support of the elements 
 of the $X_n$'s is almost disjoint, not only the support of the 
 elements  of a given sequence $(x_n)$ with 
$x_n\in X_n$, for $n\in\nat$.

  This yields the following: 
   If $X$ is a subspace of $L_p$, and
 $X$ is   not contained in $\ell_2\oplus \ell_p$,  then $X$ must
 contain a complemented copy of $\ell_p(\ell_2)$. Moreover, it 
 admits a projection onto a subspace isomorphic to $\ell_p(\ell_2)$, whose norm is arbitrarily close  to that  
 of the minimal norm projection of $L_p$ onto any subspace isomorphic to $\ell_2$.

\begin{thmB}
Let $X\subseteq L_p$ $(2<p<\infty)$. If $X$ does not embed into
$\ell_p\oplus\ell_2$ then for all $\ep>0$, $X$ contains a subspace
$Y$, which is $1+\ep$-isomorphic to $\ell_p$ $(\ell_2)$, and $Y$ is
complemented in $L_p$ by a projection of norm not exceeding $(1+\ep)
\gamma_p$ where $\gamma_p = \|x\|_p$,  $x$ being a symmetric $L_2$
normalized Gaussian random variable.

Moreover,   we  can write $Y$ as the complemented sum of\  $Y_n$'s where
 $Y_n$ is  $(1+\ep)$-isomorphic to $\ell_2$ and $Y$ is $(1+\ep)$-isomorphic
 to the $\ell_p$-sum of the $Y_n$'s, and there exists a sequence $(A_n)$ 
 of disjoint measurable sets so that $\|y|_{A_n}\|_p\ge (1-\ep2^{-n})\|y\|$
  for all $y\in Y_n$ and $n\in\nat$.
\end{thmB}

The original proof of the \cite{JO2} result about quotients of
subspaces of $\ell_p\oplus\ell_2$, is quite complicated, and a
byproduct of our results will be to give a much easier proof (see
Section~7). In addition, we can characterize when $X\subseteq L_p$
$(2<p<\infty)$ embeds into $\ell_p\oplus \ell_2$ in terms of its
asymptotic structure \cite{MMT}. {From} the \cite{KP} and \cite{JO1}
results we first note that $X\subseteq L_p$ $(2<p<\infty)$ embeds
into $\ell_p$ if and only if it is asymptotic $\ell_p$, and $X$
embeds into $\ell_2$ if and only if it is asymptotic $\ell_2$.

Let us say $X$ is {\em asymptotic $\ell_p\oplus\ell_2$} if for some
$K$ and all $(e_i)_1^n\in \{X\}_n$, the $n^{\text{th}}$ asymptotic
structure of $X$, there exists $(w_i)_1^n\subseteq [0,1]$ so that
for all $(a_i)_1^n\subseteq \real$,
\begin{equation}\label{E:1}
\Big\|\sum_1^n a_i e_i\Big\| \Ksim \Big(\sum_1^n
|a_i|^p\Big)^{1/p} \vee \Big( \sum_1^n |a_i|^2
|w_i|^2\Big)^{1/2}\ .
\end{equation}
We note that the space $\ell_p\oplus\ell_2$ is itself asymptotic $\ell_p\oplus \ell_2$.
Indeed, denote by  $(f_i)$ and~$(g_i)$ the unit vector bases of $\ell_p$ and $\ell_2$,
 respectively, viewed 
 as elements of $\ell_p\oplus \ell_2$.
  For $(x,y)\in\ell_p\oplus\ell_2$ we put $\|(x,y)\|=\|x\|_p\!\vee\!\|y\|_2$.
  Since $(f_i)$ and $(g_i)$ are $1$-subsymmetric and $\ell_p\oplus \ell_2$ is reflexive,
  the elements of the $n^{\text{th}}$ asymptotic structure of $\ell_p\oplus \ell_2$ are exactly  the sequences 
  $(z_i)_{i=1}^n$
  in $\ell_p\oplus \ell_2$, for which there are 
  $0\!=\!k_0\!<\!k_1\!<\!k_2\!<\!\ldots k_n$ in $\nat$, and  $(a_j)$, $(b_j)$ in $\real$ with
  $$z_i=\sum_{j=k_{i-1}+1}^{k_i} (a_j f_j +b_j g_j),$$
  so that $\|z_i\|=v_i\vee w_i=1$, where 
  $$v_i=\Big(\sum_{j=k_i-1}^{k_i} |a_j|^p\Big)^{1/p},\text{ and }
 w_i=\Big(\sum_{j=k_i-1}^{k_i} |b_j|^2\Big)^{1/2}.$$
 For $(\xi_i)_{i=1}^n\subset[-1,1]$ we therefore compute
 $$\Big\|\sum_{i=1}^n \xi_i z_i\Big\|=\Big(\sum_{i=1}^n |\xi_i|^p v_i^p\Big)^{1/p}
 \vee \Big(\sum_{i=1}^n |\xi_i|^2 w_i^2\Big)^{1/2}\le
 \Big(\sum_{i=1}^n |\xi_i|^p \Big)^{1/p}
 \vee \Big(\sum_{i=1}^n |\xi_i|^2 w_i^2\Big)^{1/2}.$$
 Assuming now that (otherwise \eqref{E:1} follows immediately)   $$\Big(\sum_{i=1}^n  |\xi_i|^p v_i^p \Big)^{1/p}
 \ge \Big(\sum_{i=1}^n |\xi_i|^2 w_i^2\Big)^{1/2},$$
 we  deduce that 
 $$\Big\|\sum_{i=1}^n \xi_i z_i\Big\|^p\ge \frac12\Big[\sum_{i=1}^n |\xi_i|^p v_i^p
 + \Big(\sum_{i=1}^n |\xi_i|^2 w_i^2\Big)^{p/2}\Big]\ge
 \frac12\sum_{i=1}^n |\xi_i|^p( v_i^p\vee w_i^p)=
\frac12\sum_{i=1}^n |\xi_i|^p.$$
It follows therefore that $(z_i)$ satisfies \eqref{E:1} with $K=2$ and we deduce that
$\ell_p\oplus\ell_2$ is asymptotic $\ell_p\oplus\ell_2$. 

For $n\in\nat$  let $(e_{i,j}^{(n)}:i,j\le n)$ be the unit vector basis of $\ell_p^n(\ell_2^n)$, i.e. 
$$\Big\|\sum_{i,j=1}^n a_{i,j} e^{(n)}_{i,j} \Big\|=
\Big(\sum_{i=1}^n \Big(\sum_{j=1}^n |a_{i,j}|^2\Big)^{p/2}\Big)^{1/p},
 \text{ for all $(a_{i,j})\subset \real$.}$$
 
Note that   $(e_{i,j}^{(n)})$ is, ordered lexicographically, isometrically
in the $(n^2)^{th}$ asymptotic structure of $\ell_p(\ell_2)$, for all $n\in\nat$, but it is not hard to deduce from the aforementioned description of the asymptotic structure of $\ell_p\oplus \ell_2$, that
  $(e_{i,j}^{(n)})$  is not (uniformly in $n\in\nat$) in the $(n^2)^{th}$ asymptotic structure of 
  $\ell_p\oplus\ell_2$. Theorem  B yields therefore the following  
\begin{corC}
$X\subseteq L_p$ $(2<p<\infty)$ embeds into $\ell_p\oplus\ell_2$ if
and only if $X$ is asymptotic $\ell_p\oplus \ell_2$.
\end{corC}

Indeed, if $X$ does not embed into $\ell_p\oplus\ell_2$, thenby Theorem B it
contains an isomorph of $\ell_p$ $(\ell_2)$, which is 
not  asymptotic $\ell_p\oplus\ell_2$.

Using Theorem A and Theorem B  we will be able to deduce the following
additional surprising  characterization of subspaces of $L_p$ which embed into $\ell_p\oplus\ell_2$. 
  It is analogous to the characterization of subspaces
    of $L_p$ which embed in $\ell_p$ via normalized weakly null sequences 
    (see the aforementioned result from \cite{J}) and we thank W.~B.~Johnson for having pointed it out to us.
    
    \begin{corD}
    $X\subseteq L_p$ $(2<p<\infty)$ embeds into $\ell_p\oplus\ell_2$ if and only if
     there exists a $K\ge 1$  so that
 every  normalized weakly null sequence  in $S_X$ admits a subsequence  $(x_i)$ satisfying
   for all scalars $(a_i)$, 
\begin{equation}\label{E:1.4}
\Big\|\sum a_i x_i\Big\|\Ksim\Big(\sum|a_i|^p\Big)^{1/p}\vee \Big(\sum a_i^2\|x_i\|_2^2\Big)^{1/p}.
\end{equation}

    \end{corD}
    A  proof of Corollary D will be given at the end of Section~5. It is worth noting that \eqref{E:1.4}
 is a reformulation of  \eqref{eqA1} in (b) of Theorem A. The difference here is that the constant $K$
 is uniform and not dependent on the particular sequence. Without the uniformity  assumption, the Corollary would be false
 (see Theorem \ref{thm2.4} below).
In Section~2 we recall some inequalities for unconditional basic
sequences and martingales in $L_p$. Section~3 contains the proof of
  Theorem~A, along with the necessary
preliminaries on weakly null trees, and the ``infinite asymptotic
game.''  In Section~4 we introduce a dichotomy of Kadets--Pe{\l}czynski 
type and apply the results of Section~2 to embed a class of
subspaces of $L_p$ into $\ell_p\oplus \ell_2$. Section~5 considers the
subspaces of $L_p$ which do not embed in $\ell_p\oplus \ell_2$; we
show that such subspaces contain  ``thinly supported $\ell_2$'s''.
 More precisely, for some $K<\infty$,  we find  subspaces
 $Y_n$,  $n\in\nat$, which are $K$-isomorphic to $\ell_2$, but for which the natural
 equivalence of $\|\cdot\|_p$ and $\|\cdot\|_2$ on $Y_n$ is bad. 
By this we mean that $\|y\|_p \ge M_n \|y\|_2$, for all $y\in Y_n$,
 for some sequence $(M_n)\subset\real$, with $M_n\nearrow \infty$, as $n\nearrow\infty$. This will
 enable us to argue that we can choose the $Y_n$'s so that vectors $y_n\in S_{Y_n}$, $n\in\nat$, are
 almost disjointly supported and hence the closed linear span of the $Y_n$'s is isomorphic to
 $\ell_p(\ell_2)$.
 Section~6 refines the result of Section~5, obtaining alomst
disjointly  supported $\ell_2$'s, by
applying techniques from Aldous's paper \cite{A} on random measures.
As well as the new proof of the result from \cite{JO2} mentioned
above, Section~7 includes a construction of  subspaces of
$L_p$, isomorphic to $ell_2$, which embed only with bad constants in spaces of the form
$\ell_p\oplus\big(\bigoplus_{i=1}^m \ell_2\big)_p$.  
 In Section~8 we recall what is known and not known about small 
 $\L_p$-spaces and raise  a problem about when $X\subset L_p$ embeds 
  into $\ell_p(\ell_2)$. In light of the deep work of \cite{BRS} in constructing
   uncountably many separable $\L_p$ spaces, it is likely that
   further study of their ordinal index will be needed to make progress
    on classifying the next group of  smaller $\L_p$-spaces.
 
We are especially grateful to the referee for two incredibly detailed reports
which greatly improved our exposition.


\section{Some inequalities in $L_p$}\label{S:2}    

We first recall the well known fact that an unconditional basic
sequence in $L_p$ is trapped between $\ell_p$ and $\ell_2$.

\begin{prop}\label{prop2.1}
{\em (see e.g. \cite{AO})} Let $(x_i)$ be a normalized
$\lambda$-unconditional basic sequence in $L_p$ $(2<p<\infty)$. Then
for all $(a_i)\subseteq \real$
\begin{equation*}
\lambda^{-1} \left(\sum |a_i|^p\right)^{1/p} \le \Big\|\sum a_i
x_i\Big\|_p \le \lambda\, B_p \left( \sum |a_i|^2\right)^{1/2}\ .
\end{equation*}
\end{prop}

In  Proposition \ref{prop2.1}, $B_p$ is the Khintchin constant, $\|\sum a_i
r_i\|\le B_p (\sum |a_i|^2)^{1/2}$, where $(r_i)$ is the Rademacher
sequence.

H.~Rosenthal proved that if the $x_i$'s are  independent and mean
zero random variables in $L_p$ then they span a subspace of
$\ell_p\oplus\ell_2$.

\begin{thm}\label{thm2.2}
\cite{R}  Let $2<p<\infty$. There exists $K_p <\infty$ so that if
$(x_i)$ is a normalized mean zero sequence of independent random
variables in $L_p$, then for all $(a_i)\subseteq \real$
\begin{equation*}
\Big\|\sum a_i x_i\Big\|_p  \Kpsim \left(\sum |a_i|^p\right)^{1/p}
\vee \left( \sum |a_i|^2 \|x_i\|_2^2\right)^{1/2}\ .
\end{equation*}
\end{thm}

D.~Burkholder extended this result to martingale difference
sequences as follows.

\begin{thm}\label{thm2.3}\label{B-R-I}
{\rm (\cite{B}, \cite{BDG}, \cite{Hit})} Let $2<p<\infty$. There
exists $C_p <\infty$ so that if $(z_i)$ is a martingale difference
sequence in $L_p$, with respect to the sequence $(\F_n)$ of
$\sigma$-algebras, then
\begin{equation*}
\Big\|\sum z_i\Big\|_p \Cpsim \left(\sum \|z_i\|_p^p\right)^{1/p} \vee
\Big\| \left( \sum \E [z_{i}^2|\F_{i-1}]\right)^{1/2} \Big\|_p\ ,
\end{equation*}
where $\E(x|\F)$ denotes the conditional expectation  of an integrable random variable
$x$ with respect to a sub-$\sigma$-algebra $\F$.
\end{thm}

{From} \cite{KP}, it follows that every normalized weakly null
sequence in $L_p$ admits a subsequence $(x_i)$, which is either
equivalent to the unit vector basis of $\ell_p$ or equivalent to the
unit vector basis of $\ell_2$. The latter occurs if $\ep = \lim_i
\|x_i\|_2 >0$ and the lower $\ell_2$ estimate is (essentially)
\begin{equation*}
\ep \left( \sum |a_i|^2\right)^{1/2} \le \Big\|\sum a_i x_i\Big\|_p\ .
\end{equation*}

Using Theorem~\ref{thm2.3}, W.B.~Johnson, B.~Maurey, G.~Schechtman,
and L.~Tzafriri obtained a quantitative improvement.

\begin{thm}\label{thm2.4}
\cite[Theorem 1.14]{JMST} Let $2<p<\infty$. There exists $D_p <\infty$ with the
following property. Every normalized weakly null sequence in $L_p$
admits a subsequence $(x_i)$ satisfying for some $w\in [0,1]$, for
all $(a_i)\subseteq \real$,
\begin{equation*}
\Big\|\sum a_i x_i \Big\|_p \Dpsim \left( \sum |a_i|^p \right)^{1/p}
\vee w \left(\sum |a_i|^2\right)^{1/2}\ .
\end{equation*}
\end{thm}

Thus in particular $[(x_i)]$, the closed linear subspace generated
by $(x_i)$ uniformly embeds into $\ell_p\oplus\ell_2$.

\section{A criterion for embeddability in $\ell_p\oplus\ell_2$}\label{S:3}
\bigskip

 In this section we prove Theorem A, and thus  provide an intrinsic characterization of
 subspaces of $L_p$ which  isomorphically embed into
  $\ell_p\oplus\ell_2$. This characterization is based
   on methods developed in \cite{OS1} and \cite{OS2}.

We will need the following notation.

Let $Z$ be a Banach space with a finite
 dimensional decomposition
 (FDD) $E=(E_n)$.  For $n\in\nat$, we
denote the $n$-th {\em coordinate projection } by $P^E_n$, i.e.
$P_n^E:Z\to E_n$ with  $P_n^E(z)= z_n$,
 for $z=\sum z_i\in Z$, with $z_i\in E_i$, for all $i\in\nat$.
  For a  finite
$A\subset \nat$ we put $P^E_A=\sum_{n\in A} P_n^E$.

$ c_{00}$ denotes the vector space of sequences in $ \real$ which
are eventually $0$ with unit vector basis $(e_i)$. More generally, if $(E_i)$ is a sequence of
finite dimensional Banach spaces, we define the vector space
$$ c_{00}(\oplus_{i=1}^\infty E_i)=
\Big\{(z_i):z_i\in E_i, \text{ for $i\in\nat,$ and $\{i\in\nat:
z_i\not=0\}$ is finite}\big\}.$$
The linear space        $c_{00}(\oplus_{i=1}^\infty E_i)$
  is dense in each Banach
space for which $(E_n)$ is an FDD. If $A\subset \nat$ is finite we
denote  by $\oplus_{i\in A} E_i$ the linear subspace of $
c_{00}(\oplus E_i)$ generated by the elements of $(E_i)_{i\in A}$. A
{\em blocking of $(E_i)$} is a sequence $(F_i)$ of finite
 dimensional spaces for which there is an increasing sequence
 $(N_i)$ in $\nat$ so that $(N_0=0)$
  $F_i=\oplus_{j=N_{i-1}+1}^{N_i} E_j$, for any $i\in\nat$.

Let $V$ be a Banach space with a normalized 1-unconditional basis
$(v_i)$
 and $E=(E_i)$ a sequence of finite dimensional spaces. Then
we define for $\xb=(x_i)\in c_{00}(\oplus_{i=1}^\infty E_i)$
$$\|\xb\|_{(E,V)}=\Big\|\sum_{i=1}^\infty \|x_i\|\cdot v_i\Big\|_V.$$
$\|\cdot\|_{(E,V)}$ is a norm on $ c_{00}(\oplus_{i=1}^\infty E_i)$,
and we denote the completion of   $ c_{00}(\oplus_{i=1}^\infty E_i)$,
with respect to $\|\cdot\|_{(E,V)}$, by $\big(\oplus_{i=1}^\infty
E_i\big)_V$.

For $z\in  c_{00}(\oplus E_i)$ we define the $E$-{\em support of
$z$} by  $\supp_E(z)=\{i\!\in\!\nat:  P^E_i(z)\!\not=\!0\}$. A
non-zero sequence $(z_j)\subset  c_{00}(\oplus E_i)$  is called a
{\em block sequence of $(E_i)$} if $\max\supp_E(z_n)<
\min\supp_E(z_{n+1})$, for all  $n\in\nat$,  and it is called a {\em
skipped block sequence of $(E_i)$} if $1<\min\supp_E(z_1)$ and
$\max\supp_E(z_n) < \min\supp_E(z_{n+1})-1$, for all $n\in\nat$.
Let $\deltab=(\delta_n)\subset (0,1]$. If $Z$ is a space with an FDD
$(E_i)$, we call  a  sequence $(z_j)\subset S_Z=\{z\in Z: \|z\|=1\}$
a $\deltab$-{\em skipped block sequence of $(E_n)$}, if there are
$1\!\le\! k_1\!<\! \ell_1\!< k_2\!<\!\ell_2\!<\!\cdots$ in $\nat$ so that
 $\| z_n - P^E_{(k_{n},\ell_n]}(z_n)\|\!<\!\delta_n$,
 for all $n\!\in\!\nat$.
  Of course one
could generalize the notion of  $\deltab$-skipped block sequences to
more general sequences, but we prefer to introduce this notion only
for normalized sequences. It is important to note that, in the
definition of $\deltab$-skipped block sequences, $k_1\!\ge\!1$, and, thus,
that the $E_1$-coordinate of $z_1$ is small (depending on
$\delta_1$). Let
\begin{equation*}
T_\infty=\bigcup_{\ell\in\nat}\big\{(n_1,n_2,\ldots,n_\ell)
 :n_1<n_2<\cdots n_\ell \text{ are  in $\nat$}\big\}\ .
\end{equation*}
$T_\infty$ is naturally partially ordered by extension, i.e.,
 $(m_1,m_2,\ldots m_k)\preceq(n_1,n_2,\ldots n_\ell)$ if
 $k\le \ell$ and $n_i=m_i$, for $i\le k$. 
 We call $\ell$ the length of $\alpha=(n_1,n_2,\ldots n_\ell)$
 and denote it by $|\alpha|$, with $|\emptyset|=0$
In this paper {\em trees  in a Banach space} $X$  are families in
$X$ indexed by $T_\infty$.

For a  tree $(x_\alpha)_{\alpha\in T_\infty}$ in $X$, and
$\alpha=(n_1,n_2,\ldots, n_\ell)\!\in\!T_\infty\!\cup\!\{\emptyset\}$, we
call the sequences~of the form $(x_{(\alpha,n)})_{n>n_\ell}$ {\em
nodes of $(x_\alpha)_{\alpha\in  T_\infty}$}.
 The sequences $(y_n)$,
with $y_i=x_{(n_1,n_2,\ldots,n_i)}$,  for $i\in\nat$,  for some
strictly increasing sequence $(n_i)\subset\nat$, are called {\em
branches of } $(x_\alpha)_{\alpha\in T_\infty}$. Thus, branches of a
tree $(x_\alpha)_{\alpha\in T_\infty}$ are sequences of the form
$(x_{\alpha_n})$ where $(\alpha_n)$ is a maximal linearly ordered
(with respect to extension) subset of $T_\infty$.

If $(x_\alpha)_{\alpha\in T_\infty}$  is a tree in $X$ and if
$T'\subset T_\infty$  is closed under taking initial segments
(if $(n_1,n_2,\ldots, n_\ell)\in T'$ and $m<\ell$ then 
 $(n_1,n_2,\ldots, n_m)\in T'$)
 and has
the property that for each $\alpha\!\in\!T'\cup\{\emptyset\}$
infinitely many direct successors of $\alpha$ are also in $T'$ then
we call $(x_\alpha)_{\alpha\in T'}$ a {\em full subtree of}
$(x_\alpha)_{\alpha\in T_\infty}$. Note that $(x_\alpha)_{\alpha\in
T'}$ could then be relabeled to a family indexed by $T_\infty$ and
note that the branches of $(x_\alpha)_{\alpha\in T'}$ are branches
of $(x_\alpha)_{\alpha\in T_\infty}$ and that the nodes of
$(x_\alpha)_{\alpha\in T'}$ are subsequences of  certain nodes of
$(x_\alpha)_{\alpha\in T_\infty}$.

We call a tree $(x_\alpha)_{\alpha\in T_\infty}$  in $X$ {\em
normalized} if $\|x_\alpha\|\!=\!1$, for all  $\alpha\in T_\infty$ and
{\em weakly null} if every node is a weakly null sequence. If $X$
has an $FDD$ $(E_i)$ we call
 $(x_\alpha)_{\alpha\in T_\infty}$ a {\em block tree with respect to $(E_i)$ } if every
  node and every branch $(y_n)$ is a block sequence with respect to $(E_i)$.

Note that, if $(E_i)$ is an  FDD for $X$ and if
$(\ep_\alpha)_{\alpha\in T_\infty}\!\subset\!(0,1)$, every normalized
weakly null tree $(x_\alpha)_{\alpha\in T_\infty}\!\subset\!X$
 has a full subtree  $(z_\alpha)_{\alpha\in T_\infty}$ which is an $(\ep_\alpha)$-perturbation
 of a block tree $(y_\alpha)$  with respect to $(E_i)$, i.e.
 $\|z_\alpha\!-\!y_\alpha\|\!\le\!\ep_\alpha$, for any $\alpha\!\in\!T_\infty$.
  Let us also mention that the proof of the fact, that normalized weakly null sequences
   have basic subsequences whose basis constants are arbitrarily close to 1,  generalizes
    to trees. This means that for a given $\ep\!>\!0$, and for any 
     Banach space $X$, every normalized weakly null tree in $X$ has a full subtree,
     all of whose nodes and all of whose branches are basic, and their basis constant  
      does not exceed~$1\!+\!\ep$.

 We now can state the main results of this section.
\begin{thm}\label{T:3.1}
Let $X$ be a subspace   of $L_p$, $2<p<\infty$, and assume that
there is a $C>1$ so that every normalized weakly null tree in $X$
admits a branch $(y_i)$ for which
$$\Big\|\sum_{i=1}^\infty a_i y_i\Big\|_p\Csim
\max\Bigg(\Big(\sum_{i=1}^\infty |a_i|^p\Big)^{1/p},
\Big\|\sum_{i=1}^\infty a_i y_i\Big\|_2\Bigg)\text{  for all 
$(a_i)\in c_{00}$}.$$
 Then there
is a blocking $H=(H_n)$ of the Haar basis $(h_n)$  so that
$$T: X\to \ell_p\oplus L_2,\quad  T(x)= \big((P^H_n(x))_{n\in\nat},
x\big)\in \big(\oplus_{n=1}^\infty H_n\big)_{\ell_p}\oplus
L_2\hookrightarrow
 \ell_p\oplus L_2,$$
 is an isomorphic embedding. 
\end{thm}
Theorem \ref{T:3.1} is a special case of the following result. By a $1$-subsymmetric basis we mean one that is
 $1$-unconditional and $1$-spreading.
\begin{thm}\label{T:3.2}
Let $X$ and $Y$  be separable Banach spaces, with $X$  reflexive. Let $V$ be
 a Banach space with a $1$-subsymmetric  and normalized basis $(v_i)$, and
 let $T:X\to Y$ be linear and bounded.

Assume that  for some $C\ge1$ every normalized  weakly null tree of
$X$ admits a branch  $(x_n)$
 so that
\begin{equation}\label{E:3.2.1}
\Big\| \sum_{i=1}^\infty a_n x_n\Big\|_X\Csim \Big\|
\sum_{i=1}^\infty a_n v_n\Big\|_V \vee \Big\|T\Big(\sum_{i=1}^\infty
a_n x_n\Big)\Big\|_Y \text{ for all $(a_i)\in c_{00}$}.
\end{equation}

Then there is a sequence of finite dimensional spaces $(G_i)$, so
that $X$
 is isomorphic to a subspace of $\big(\oplus_{i=1}^\infty G_i\big)_V\oplus Y$\!\!.

More precisely, under the above assumptions, if $Z$ is any reflexive space with an FDD $(E_i)$,
 and if
 $S: X\to Z$ is an isomorphic embedding, then there is a blocking
$(G_i)$ of $(E_i)$ so that $S$ is a bounded linear operator from $X$
to $\big(\oplus_{i=1}^\infty G_i\big)_V$ and the operator
$$(S,T): X\to  \big(\oplus_{i=1}^\infty G_i\big)_V\oplus Y,\quad
x\mapsto \big(S(x),T(x)\big),$$ is an isomorphic embedding.
\end{thm}

\begin{remark} Theorem \ref{T:3.1} can be obtained from Theorem  \ref{T:3.2}
  by letting $V=\ell_p$,  $Y=L_2$, $Z=L_p$, with the FDD $(E_i)$  given by the Haar
  basis,  $S$ is the inclusion map from $X$ into $L_p$ and $T$  is the formal identity map from $L_p$ to $L_2$
  restricted to $X$.

As noted in \cite[Corollary 2, Section 2]{OS2} (see also \cite{OS1}
for similar versions) the tree condition in Theorem \ref{T:3.2} can
be interpreted as follows in terms of the ``infinite asymptotic
game'', (IAG) as it has been called by Rosendal \cite{Ro}.

Let $C\ge1$ and let $\cA^{(C)}$ be the set of all sequences $(x_n)$
in $S_X$ which are $C$-basic
 and satisfy condition  \eqref{E:3.2.1}. The (IAG) is played by two players:
Player I chooses a subspace $X_1$ of $X$ having  finite
co-dimension, and  Player II chooses $x_1\in S_{X_1}$, then, again
Player I chooses a subspace $X_2$ of $X$ of finite codimension , and
Player II chooses an $x_2\in S_{X_2}$. These moves are repeated
infinitely many times,
 and Player I is declared the winner of the game if the resulting sequence $(x_n)$ is in $\cA^{(C)}$.

 $\cA^{(C)}$ is closed with respect to the infinite product of $(S_X,d)$,
 where $d$ denotes the discrete topology on $S_X$. This implies that 
 this game is determined \cite{Ma}, i.e.,
either Player I or Player II has a winning strategy and as noticed
in
   \cite[Corollary 2, Section 2]{OS2} for all $\ep>0$  Player I has a winning strategy
 for $\cA^{(C+\ep)}$ if and only  if for all  $\ep>0$, every weakly null tree in $S_X$ has a branch,
 which lies in  $\cA^{(C+\ep)}$.
\end{remark}

\begin{proof}[Proof of Theorem A  using Theorem \ref{T:3.1}]
The interpretation of our tree condition in terms of the infinite
asymptotic  game, easily
 implies  that the existence of a uniform
$C\ge 1$, so that all weakly null trees $(x_\alpha)\subset S_X$
admit a branch in $A^{(C)}$, is equivalent to the condition, that
every weakly null tree $(x_\alpha)\subset S_X$   admits a branch in
$\cA^{(C)}$, for some $C\ge 1$.

Indeed, if such a uniform $C$ does not exist, Player II could choose
a sequence $(C_n)$ in $ \real^+$ which increases to $\infty$ and
  could play the following strategy: first he follows his  winning strategy for achieving a sequence
 $(x_n)$ outside of $\cA^{(C_1)}$ and after
 finitely many steps, $s_1$, he  must have chosen a sequence
$x_1, x_2,\ldots, x_{s_2}$, which is either not $C_1$-basic or  does
not satisfy   \eqref{E:3.2.1}
 for some $a=(a_i)_{i=1}^{s_1}\in \real^{s_1}$. Then Player II follows his strategy for getting a sequence
 outside of $\cA^{(C_2)}$, and continues that way using $C_3$, $C_4$ etc.
 It follows that the infinite sequence $(x_n)$, which is obtained
by Player II cannot be in any $\cA^{(C)}$. Therefore Player II  has
a winning strategy for choosing
 a sequence outside of $\bigcup_{C\ge1}\cA^{(C)}$ which
 means that there is a weakly null tree, $(z_\alpha)$, none of whose
branches are in  $\bigcup_{C\ge1}\cA^{(C)}$ .

Using Theorem \ref{T:3.1}, we deduce therefore (b)$\Rightarrow$(a) of
Theorem A. The implication (a)$\Rightarrow$(c) in  Theorem A is
 easy, using arguments like those above establishing that $\ell_p\!\oplus\!\ell_2$ is
  asymptotic $\ell_p\oplus\ell_2$.

In order to show
 (c)$\Rightarrow(b)$ let $(x_\alpha)$ be a normalized weakly null tree in $L_p$.
 After passing to a full subtree, and perturbing, we can assume that
  $(x_\alpha)$ is a block tree with respect to the Haar basis.
 By (c) there is branch $(z_n)$, a sequence $(w_i)\subset[0,1]$ and
 $C\ge 1$
  so that
 \begin{align}\label{E:3.1a}
\Big\|\sum a_i z_i\Big\|_p\Csim \Big(\sum |a_i|^p\Big)^{1/p}\vee
  \Big(\sum w_i^2a_i^2\Big)^{1/2}\text{ for all $(a_i)\in c_{00}$}.
\end{align}
Since $(z_i)$ is an unconditional sequence  and since
$\|\cdot\|_2\le \|\cdot\|_p$ on $L_p$ it follows from
 Proposition \ref{prop2.1}
that for some constant $c_p$
\begin{align}\label{E:3.1}
\Big\|\sum a_i z_i\Big\|_p\ge c_p\Big(\sum |a_i|^p\Big)^{1/p}\vee
\Big\|\sum a_i z_i\Big\|_2.
\end{align}
We claim that   our branch $(z_n)$  satisfies \eqref{eqA1} for some $K<\infty$.
 Assuming this were not true, then  we could use \eqref{E:3.1a},  and  choose a normalized block sequence $(y_n)$
 of $(z_n)$, say
$$y_n=\sum_{i=k_{n-1}+1}^{k_n} a_i z_i, \text{ with $a_i\in\real$, for $i\in\nat$ and
$0=k_0<k_1<\ldots$,}$$ so that for all $n\in\nat$
\begin{align}\label{E:3.2}
&\sum_{i=k_{n-1}+1}^{k_n} w_i^2a_i^2=1,\text{ and}\\
\label{E:3.3}&\Big(\sum_{i=k_{n-1}+1}^{k_n} |a_i|^p\Big)^{1/p}\vee
\|y_n\|_2< 2^{-n}.
\end{align}
For any $(b_i)\in c_{00}$  it follows therefore from \eqref{E:3.1a}
that
$$
\Big\|\sum b_n y_n\Big\|_p\Csim\Big(\sum |b_n|^2     \Big)^{1/2},
$$
thus $(y_n)$ is $C$-equivalent to the unit vector basis of $\ell_2$.
The result by
Kadets and Pe{\l}czy\'nski \cite{KP} yields that $\|\cdot\|_p$ and $\|\cdot\|_2$
 must be equivalent on $Y$. But $\lim_{n\to\infty}\|y_n\|_2=0$ by \eqref{E:3.3}, so 
 we have a contradiction.
\end{proof}

For the proof of Theorem \ref{T:3.2} we need to recall some
 results from \cite{OS1} and \cite{OS2}.
The following result restates Corollary 2.9 of \cite{OS2}, versions
of which where already shown in \cite{OS1}.

\begin{thm}\label{T:3.3}{\rm\cite[Corollary 2.9 (c)$\iff$(d), and
``Moreover''-part]{OS2}}\\
Let $X$ be a subspace of a reflexive space $Z$ with an $FDD$ $(E_i)$
and let
$$\cA\subset \{ (x_n): x_n\in S_X\text{ for $n\in\nat$ }\}.$$
Then the following are equivalent.
\begin{enumerate}
\item[a)] For any $\epb=(\ep_n)\subset(0,1)$ every weakly null tree
 in $S_X$ admits a branch in $\overline{\cA_{\epb}}$, where
 $$\cA_{\epb}=\big\{(x_n)\subset S_X: \exists (z_n)\!\in\!\cA\quad
  \|z_n-x_n\|\le \ep_n\text{ for $n\in\nat$}\big\},$$
and where $\overline{\cA_{\epb}}$ denotes the closure in the product
 of the discrete topology on $S_X$.
 \item[b)] For any $\epb=(\ep_n)\subset(0,1)$ there is a blocking
  $(F_i)$ of $(E_i)$ so that every $c\epb$-skipped block sequence
  $(x_n)\subset S_X$  of $(F_i)$ lies in $\overline{\cA_{\epb}}.$
  Here $c\in(0,1)$ is a constant which only depends on the projection constant
  of $(E_i)$ in $Z$.
\end{enumerate}
\end{thm}

We also need a blocking lemma which appears in various forms in
\cite{KOS}, \cite{OS1}, \cite{OS2} \cite{OSZ} and ultimately results
from a blocking trick of W.~B.~Johnson \cite{J}.
 In the statement of Lemma \ref{L:3.4} (and elsewhere) reference
 is made to the {\em weak$^*$-topology of $Z$}, a space with a boundedly complete FDD
 $(E_i)$. By this we mean the weak$^*$-topology on $Z$ obtained by regarding 
 it as the dual space of the norm closure of the span of  $(E_i^*)$ in $Z^*.$
 This is then just the topology of coordinatewise convergence in $Z$ with respect to
 the coordinates of $(E_i)$.
\begin{lem}\label{L:3.4}{\rm \cite[Lemma 3, Section 3]{OS2}}
Let $X$ be a  subspace of a space $Z$ having a boundedly complete
FDD $E=(E_i)$ with projection constant $K$ with $B_X$ being a
$w^*$-closed subset of $Z$. Let $\delta_i\downarrow 0$. Then there
exist $0= N_0 < N_1<\cdots$ in $\nat$  with the
following properties. For all $x\in S_X$ there exists
$(x_i)_{i=1}^\infty \subseteq X$, and for all $i\in \nat$, there
exists $t_i\in (N_{i-1},N_i)$ satisfying ($t_0=0$ and $t_1>1$)
\begin{itemize}
\item[a)] $x = \sum_{j=1}^\infty x_j$,
\item[b)] $\|x_i\| < \delta_i$ or $\|P_{(t_{i-1},t_i)}^E x_i - x_i\|
<\delta_i \|x_i\|$,
\item[c)] $\|P_{(t_{i-1},t_i)}^E x-x_i\| < \delta_i$,
\item[d)] $\|x_i\| < K+1$,
\item[e)] $\|P_{t_i}^E x\| < \delta_i$.
\end{itemize}
\end{lem}
\begin{proof}[Proof of Theorem \ref{T:3.2}]
Assume $X$ embeds in a reflexive space $Z$ with an FDD $E=(E_i)$.
By Zippin's theorem \cite{Z} such a space $Z$ always exists. After renorming we can assume
that the projection constant $K=\sup_{m\le n}\|P^E_{[m,n]}\|=1$ and that $X$ is (isometrically)
 a subspace of $Z$.
  We
also assume without
 loss of generality that $\|T\|= 1$.

For a  sequence $\xb=(x_i)\in S_X$ and $a=\sum a_i e_i\in c_{00}$ we
define
$$\Btrivert\sum a_i e_i\Btrivert_{\xb}=
\Big\|\sum a_i v_i\Big\|_V\vee \Big\|T\Big(\sum a_i
x_i\Big)\Big\|_Y.$$ Then $\trivert\cdot\trivert_{\xb}$ is a norm on $
c_{00}$ and
 we denote the completion of $ c_{00}$ with respect to  $\trivert\cdot\trivert_{\xb}$ by $W_{\xb}$.

Define
$$\cA=\left\{\xb=(x_n)\subset S_X: \begin{matrix} \xb \text{ is $\frac32$-basic and $\frac32 C$-equivalent}\\
                                                  \text{to $(e_i)$ in $W_{\xb}$}
                                   \end{matrix}\right\}.$$

Observe that
condition a) of Theorem \ref{T:3.3} is satisfied for this set $\cA$. Indeed, given any
weakly null tree in    
$S_X$
 we may assume, as noted before the statement of
Theorem \ref{T:3.1} that, by passing to a full subtree,  the branches are basic
with a constant close to $1$, and, thus the first requirement of the definition
of $\cA$ can be satisfied. The hypothesis from Theorem \ref{T:3.2}
then guarantees
that
$\cA_{\overline{\ep}}$
 contains the required branch.

     We first choose a null sequence $\epb=(\ep_i)\subset(0,1)$, which decreases fast enough                
    to $0$ to ensure that every sequence $\xb=(x_n)$ in $\overline{\cA_{\epb}}$       
    is $2$-basic and $2C$ equivalent to $(e_i)$ in $W_{\xb}$.          
                   By Theorem \ref{T:3.3} applied to $\epb$ we
 can find a blocking  $F=(F_i)$ of $(E_i)$ and a sequence,
  so that every  $c\epb$-skipped block sequence $(x_i)\subset S_X$  of $(F_i)$ 
   ($c$ is the constant in Theorem \ref{T:3.3} (b))
is $2$-basic and $2C$-equivalent
 to $(e_i)$ in $W_{\xb}$. 
  We put $\deltab=(\delta_i)=c\epb$.
 Then we apply Lemma \ref{L:3.4} to get a further blocking $(G_i)$, $G_i=\oplus_{j=N_{i-1}+1}^{N_i} F_j$,
 for $i\in\nat$ and some sequence $0=N_0<N_1<N_2\ldots $, so that for every $x\in S_X$ there is a sequence $(t_i)\subset N$, with
 $t_i\in(N_{i-1},N_i)$ for $i\in\nat$,  and  $t_0=0$, and a sequence $(x_i)$ satisfying  (a)-(e).

We also may assume that $\sum_{i=1}^\infty \delta_i<1/36 C$ and
 will show that for every $x\in X$
\begin{equation}\label{E:3.2.2}
\|x\|_X\Csima{36}\Big(\Big\|\sum_{i=1}^\infty \|P_i^G(x)\|
v_i\Big\|_V\Big) \vee \|T(x)\|_Y.
\end{equation}
This implies that the map $ X\to (\oplus G_i)_V\oplus Y,\quad
x\mapsto ((P_i^G(x)),T(x))$, is an isomorphic embedding.

Let $x\in S_X$ and choose $(t_i)\subset \nat$ and $(x_i)\subset X$
as prescribed in Lemma \ref{L:3.4}. Letting
 $B=\big\{i\ge 2: \|P^F_{(t_{i-1},t_i)}(x_i)-x_i\|\le \delta_i\|x_i\|\big\}$
it follows that $(x_i/\|x_i\|)_{i\in B}$ is a $\deltab$-skipped
block sequence of $(F_i)$ and therefore

\begin{equation}\label{E:3.2.2a}
\Big\|\sum_{i\in B} x_i\Big\|_X\Csima{2}   \Big\|\sum_{i\in B}
\|x_i\| v_i\Big\|_V \vee \Big\|T\Big(\sum_{i\in B} x_i\Big)\Big\|.
\end{equation}
We want to estimate $\big\|\sum_{i=1}^\infty \|x_i\| v_i\big\|_V\vee\|T(x)\|$. Since $1\not\in B$  
(no matter how large $\|x_1\|$ is) we will distinguish between the case that $\|x_1\|$ 
is essential  and the case that $\|x_1\|$ is small enough to be discarded.

If $\|x_1\|\ge 1/8C$ then we deduce
 that
\begin{align}\label{E:3.2.3}
\frac1{8C}&\le \|x_1\|
 \le \Big\|\sum_{i=1}^\infty \|x_i\| v_i\Big\|_V\vee \|T(x)\|_Y \\
 &\le \Big(\Big\|\sum_{i\in B}^\infty \|x_i\|v_i\Big\|_V + \|x_1\|+ \sum_{i\not\in B} \delta_i\Big) \vee  \|T(x)\|_Y\notag\\
 &\le 2C\Big\|\sum_{i\in B}^\infty x_i\Big\|+2+\sum\delta_i\quad
 \text{[by \eqref{E:3.2.2a}, (d) of Lemma \ref{L:3.4}] and since $\|T\|=1$]}\notag\\
 &\le 2C\|x\|+ 2C\Big\|\sum_{i\not\in B}^\infty x_i\Big\|+2+\sum\delta_i\notag\\
 &\le 2C\|x\|+2C\|x_1\|+2C\sum\delta_i +2+\sum\delta_i \le 9C.\notag
\end{align}

If $\|x_1\|< 1/8C$ then
\begin{align*}
1&=\|x\|\le \Big\|\sum_{i\in B} x_i\Big\| +\frac1{4C}\\
 &\le 2C\Big(\Big\|\sum_{i\in B}\|x_i\| v_i\Big\|_V\vee\Big\|T\Big(\sum_{i\in B} x_i \Big)\Big\|_Y\Big)+\frac1{4C}\qquad  \text{[By \eqref{E:3.2.2a}]} \\
&\le 2C\Big(\Big\|\sum_{i=1}^\infty\|x_i\|
v_i\Big\|_V\vee\|T(x)\|_Y\Big)+ \frac12 +\frac1{4C}\le
 2C\Big(\Big\|\sum_{i=1}^\infty\|x_i\| v_i\Big\|_V\vee\|T(x)\|_Y\Big)+ \frac34
\end{align*}
and, thus,
\begin{align}\label{E:3.2.4}
\frac1{8C}&
\le  \Big\|\sum_{i=1}^\infty \|x_i\| v_i\Big\|_V\vee \Big\|T(x)\Big\|_Y \\
&\le \Big( \Big\|\sum_{i\in B} \|x_i\| v_i\Big\|_V \vee
\Big\|T\Big(\sum_{i\in B} x_i\Big)\Big\|_Y\Big) + \frac1{4C}\notag\\
&\le 2C\Big\|\sum_{i\in B} x_i\Big\|+ \frac1{4C}\quad
\text{[By \eqref{E:3.2.2a}]}\notag\\
&\le 2C\|x\|+ 2C\|x_1\|+2C\sum \delta_i +\frac1{4C}\le 8C.\notag
\end{align}
\eqref{E:3.2.3} and \eqref{E:3.2.4} imply that
\begin{equation}\label{E:3.2.5}
 1\Csima{9} \Big\|\sum_{i=1}^\infty \|x_i\| v_i\Big\|_V\vee \big\|T(x)\big\|.
\end{equation}
For $n\in\nat$ define $y_n=P^F_{(t_{n-1},t_n]}(x)$. From Lemma
\ref{L:3.4} (c) and (e) it follows that
$\|y_n-x_n\|\le\|P^F_{(t_{n-1},t_n)}(x)-x_n\|+\|P^F_{t_n}(x)\|\le
2\delta_n$ and
 thus $\sum\|y_n-x_n\|\le 1/18C$ which implies by \eqref{E:3.2.5} that
 \begin{equation}\label{E:3.2.6}
1\Csima{18} \Big\|\sum_{i=1}^\infty \|y_i\| v_i\Big\|_V\vee
\big\|T(x)\big\|.
\end{equation}
  Since for $n\in\nat$ we have $(N_{n-1},N_n]\subset(t_{n-1},t_{n+1})$ and
 and $(t_{n-1},t_n]\subset (N_{n-2},N_{n})$ (put $N_{-1}=N_0=0$ and $P_0^G=0$) it follows
from  the assumed $1$-subsymmetry of $(v_n)$ and the assumed
bimonotonicity of $(E_i)$ in $Z$ that
\begin{align*}
\frac12\Big\|\sum_{n\in\nat} \|y_n\| v_n\Big\|_V&\le
 \frac12\Big\|\sum_{n\in\nat} \big(\|P^G_{n-1}(x)\|+\|P^G_{n}(x)\|)v_n\Big\|_V\\
&\le  \Big\|\sum_{n\in\nat} \|P^G_{n}(x)\|v_n\Big\|_V\\
&\le  \Big\|\sum_{n\in\nat}  \big\|P^F_{(t_{n-1},t_{n+1})}(x)\big\|v_n\Big\|_V\\
&\le  \Big\|\sum_{n\in\nat}
\big(\|y_n\|+\|y_{n+1}\|\big)v_n\Big\|_V\le 2\Big\|\sum_{n\in\nat}
\|y_n\| v_n\Big\|_V,
\end{align*}
which implies with \eqref{E:3.2.6}  that
\begin{equation*}
1\Csima{36} \Big\|\sum_{i=1}^\infty \|P_i^G(x)\| v_i\Big\|_V\vee
\big\|T(x)\big\|.
\end{equation*}
and finishes the proof of our claim.
\end{proof}

\section{Embedding small subspaces in $\ell_p\oplus\ell_2$}\label{S:4}

For a subspace $X$ of $L_p$ (where $p>2$, as everywhere in this
paper) we shall say that a function $v$ in $L_{p/2}$ is a {\it
limiting conditional variance} associated with $X$ if there is a
weakly null sequence $(x_n)$ in $X$ such that $x_n^2$ converges to
$v$ in the weak topology of $L_{p/2}$.  It is equivalent to say
that, for all $E\in \Sigma$ (recall that $L_p$ was defined over the atomless  and separable
 probability space $(\Omega,\Sigma,\mathbb P)$)
$$
\mathbb E[1_Ex_n^2]\to \mathbb E[1_Ev]
$$
as $n\to \infty$.  The set of all such $v$ will be denoted $V(X)$.
Note that, because $p>2$, every weakly null sequence  $(x_n)$
in $X$ does of course have a subsequence $(x_{n_k})$ such that
$x_{n_k}^2$ converges (to some $v\in V(X)$) for the weak topology of
the reflexive space $L_{p/2}$.

Limiting conditional variances occur naturally in the context of the
martingale inequalities to be used in this section, and are closely
related to the random measures of Section \ref{S:6}.  It
is therefore natural to express the basic dichotomy underlying our
main Theorem  B in terms of $V(X)$.

\begin{prop}\label{K-P}
Let $X$ be a subspace of $L_p$, where $p>2$.  One of the following
is true:

{\rm (A)} there is a constant $M>0$ such that $\|v\|_{p/2}\le
M\|v\|_1$ for all $v\in V(X)$;

{\rm (B)} no such constant $M$ exists, in which case there exist
disjoint sets $A_i\in \Sigma$ and elements $v_i\in V(X)$ ($i\in
\mathbb N$), such that $\|1_{A_i}v_i\|_{p/2}\to 1$ and
$\|1_{\Omega\setminus A_i}v_i\|_{p/2}\to 0$ as $i\to \infty$.
\end{prop}
\begin{proof}
This is a  consequence of the Kadets--Pe{\l}czynski
dichotomy. Either there exists an $\ep>0$ so that
$$
V(X)\subset  \big\{u\in L_{p/2}: \mathbb P[|u|\ge \ep\|u\|_{p/2}] \ge
\ep\big\}
$$
then  
$$\|u\|_1\ge \E\big[\ep \|u\|_{p/2} 1_{[|u|\ge \ep\|u\|_{p/2}]}\big]\ge \ep^2 \|u\|_{p/2},\text{ for all $u\in V(X)$},$$
 and (A)
holds for $M=\ep^{-2}$.  Otherwise, by the construction in Theorem~2 of \cite{KP}, we
obtain (B).
\end{proof}

The rest of this section will be devoted to showing that if (A)
holds then $X$ embeds in $\ell_p\oplus \ell_2$. By Theorem~3.1, it
will be enough to prove the following proposition.

\begin{prop}\label{P:4.2}
Let $X$ be a subspace of $L_p$, where $p>2$, and assume that (A)
holds in Proposition~\ref{K-P}. Then there is a constant $K$ such that
every weakly null tree in $S_X$ has a branch $(x_i)$ satisfying
$$
K^{-1}\Big\|\sum c_ix_i\Big\|_p\le \max\Big\{\Big(\sum|c_i|^p\Big)^{1/p}, \Big\|\sum
c_ix_i\Big\|_2\Big\}\le K \Big\|\sum c_ix_i\Big\|_p,
$$
for all $c_i\in \mathbb R$.
\end{prop}
\begin{proof}
Our proof, using Burkholder's martingale version of Rosenthal's
Inequality (Theorem \ref{B-R-I}), is closely modeled on Theorem 1.14 of
\cite{JMST}. Let $(x_\alpha)_{\alpha\in T_\infty}$ be a weakly null
tree in $S_X$. Taking small perturbations, we may suppose that we
are dealing with a block tree of the Haar basis.  So for each
$\alpha\in T_\infty$, $x_\alpha$ is a finite linear combination of
Haar functions, say $x_\alpha\in [h_n]_{n\le n(\alpha)}$, and for
each successor $(\alpha,k)$ of $\alpha$ in $T_\infty$,
$x_{(\alpha,k)}\in [h_n]_{n(\alpha)<n\le n(\alpha,k)}$. We may then
proceed to choose a full subtree $T'$
 of $T_\infty$ 
  having the properties (1) and (2), below, as we now describe.

First, we consider the first level of the tree, that is to say the
sequence of elements $x_{(n)}$ with $n\in \mathbb N$.  We may
extract a subsequence for which $x_{(n)}^2$ converges weakly in
$L_{p/2}$ to some $v_0\in V(X)$ and then, by leaving out a finite
number of terms, ensure that $|\mathbb E[x_{(n)}^2]^{1/2}-\mathbb
E[v_0]^{1/2}|<\frac12$.

We now continue by taking subsequences of the successors of each
$\alpha$ in such a way that the following hold
 (for $n\in\nat$, $\mathcal H_n$ denotes the $\sigma$-algebra generated by 
  $(h_i:i\!\le\!n)$) :
 \begin{enumerate}
 \item the elements $x^2_{(\alpha,n)}$ (with $(\alpha,n)\in T'$)
  of $L_{p/2}$  converge weakly to some $v_\alpha\in V(X)$;
 \item for all $(\alpha,k)\in T'$ we have $\|\mathbb
 E[x_{(\alpha,k)}^2\mid \mathcal H_{n(\alpha)}]^{1/2}- \mathbb E[v_\alpha\mid \mathcal
 H_{n(\alpha)}]^{1/2}\|_\infty < 2^{-|\alpha|-1}$.
 \end{enumerate}
To achieve the above, we use  our earlier
remark based on relexivity of $L_{p/2}$, and  the fact that
weak convergence implies norm convergence in the finite dimensional
space $[h_n]_{n\le n(\alpha)}$.

We now take any branch $(x_i)$ of the resulting subtree
$(x_\alpha)_{\alpha\in T'}$.  So $x_i= x_{\alpha_i}$ where
$\alpha_i$ is the initial segment  $(n_1,n_2,\dots,n_i)$ of some branch
$(n_1,n_2,\dots)$ of $T'$.  We consider the $\sigma$-algebras
$\mathcal F_i$ where $\mathcal F_0=\{\emptyset, \Omega\}$ and
$\mathcal F_i = \mathcal H_{n(\alpha_i)}$ for $i\ge 1$ and write
$\mathbb E_i$ for the conditional expectation relative to $\mathcal
F_i$. Since we are dealing with a block tree the sequence $(x_i)$ is
a block basis of the Haar basis, and hence a
martingale-difference sequence with respect to $(\mathcal F_i)$. We
may therefore apply Theorem~\ref{B-R-I} to conclude that
the $L_p$-norm of a linear combination $\sum c_i x_j$ is
$C_p$-equivalent to
$$
\max\Big\{\Big(\sum |c_i|^p\Big)^{1/p}, \Big\|\sum c_i^2 \mathbb
E_{i-1}[x_i^2]\Big\|_{p/2}^{1/2}\Big\}.
$$

 We shall show
that, provided we modify the constant of equivalence, we may replace
the second term in this expression by
$$
\Big\|\sum c_i^2 \mathbb E_{i-1}[x_i^2]\Big\|_1^{1/2},
$$
which equals $\|\sum c_ix_i\|_2$.

Now, by construction, the conditional expectations $\E_{i-1}[x_i^2]$
are close to $\E_{i-1}[v_{i-1}]$, where, for $j\ge 1$, $v_j$ denotes
$v_{\alpha_j}$.  More precisely, we may use (2) above and the triangle inequality
 in $L_{p}(\ell_2)$ to obtain
\begin{equation}
\left|\Big\|\sum c_i^2\mathbb E_{i-1}[x_i^2]\Big\|_{p/2}^{1/2}\!-\!\Big\|\sum c_i^2
\mathbb E_{i-1}[v_{i-1}]\Big\|_{p/2}^{1/2}\right|
 \!\le\! \Big\|\big(\sum c_i^2 2^{-2i}\big)^{1/2}\Big\|_p\! \le\! \max |c_i|.\label{LpI}
\end{equation}
 We similarly get
\begin{equation}
\left|\Big\|\sum c_i^2\mathbb E_{i-1}[x_i^2]\Big\|_1^{1/2}-\Big\|\sum c_i^2
\mathbb E_{i-1}[v_{i-1}]\Big\|_1^{1/2}\right| \le\Big \|(\sum c_i^2
2^{-2i})^{1/2}\Big\|_2 \le \max |c_i|.\label{L2I}
 \end{equation}
 Using our assumption about
 $V(X)$,  the fact that all the $v_i$ are non-negative and the inequalities
   \eqref{LpI} and \eqref{L2I} we obtain
 \begin{align*}
 \Big\|\sum c_i^2\mathbb E_{i-1}[x_i^2]\Big\|_{p/2}^{1/2}&\le 
\Big\|\sum c_i^2 \mathbb E_{i-1}[v_{i-1}]\Big\|_{p/2}^{1/2}+\max|c_i|\\
&\le \Big(\sum c_i^2\|\E_{i-1}[v_{i-1}]]\|_{p/2}\Big)^{1/2}+\max|c_i|\\
 &\le  \Big(\sum c_i^2\big\| v_{i-1}\big\|_{p/2}\Big)^{1/2}+\max|c_i|\\
 &\le \sqrt{M}\Big(\sum c_i^2\big\|v_{i-1}\big\|_{1}\Big)^{1/2}+\max|c_i|\\
 &=\sqrt{M}\Big\|\sum c_i^2 \mathbb E_{i-1}[v_{i-1}]\Big\|_{1}^{1/2}+\max|c_i|\\
 &\le\sqrt{M} \Big\|\sum c_i^2\mathbb E_{i-1}[x_i^2]\Big\|_1^{1/2}+\big(1+\sqrt{M}\big)\max|c_i|
 \end{align*}
 which yields the left most inequality in Proposition  \ref{P:4.2}. The right hand inequality
  is easy by Proposition \ref{prop2.1}  since $\|\cdot\|_p\ge \|\cdot\|_2$ and $(x_i)$ is unconditional, being  a block basis of the Haar basis. 
\end{proof}

\begin{cor}
Let $X$ be a subspace of $L_p$, where $p>2$, and assume that (A)
holds in Proposition~\ref{K-P}.  Then $X$ embeds isomorphically into
$\ell_p\oplus\ell_2$.
\end{cor}


\section{Embedding $\ell_p(\ell_2)$ in $X$}\label{S:5}

\begin{thm}\label{SkinnyTed}
Let $X$ be a subspace of $L_p$ $(p>2)$ and suppose that (B) of
Proposition~\ref{K-P} holds.   Then $X$ contains a subspace isomorphic to $\ell_p(\ell_2)$.
\end{thm}

The first step in the proof is to find $\ell_2$-subspaces of $X$
which have ``thin support''.  The precise formulation of this notion
that we shall use in the present section is given in the following
lemma.

\begin{lem}\label{Thinell2}
Suppose that (B) of Proposition~\ref{K-P} holds.  Then, for every $M>0$ there is
an infinite-dimensional subspace $Y$ of $X$, on which the $L_p$ and
$L_2$ norms are equivalent, but in such a way that $\|y\|_p\ge
M\|y\|_2$ for all $y\in Y$.
\end{lem}
\begin{proof}
By hypothesis, for every $M'>0$ there exists $v\in V(X)$ such that
$\|v\|_1=1$ and $\|v\|_{p/2}>M'^2$.  There is a weakly null sequence
$(x_n)$ in $X$ such that $x_n^2$ converges weakly to $v$ in
$L_{p/2}$.  By taking small perturbations of the $x_n$'s (with respect
 to the $L_p$-norm) and  by noting that  the Cauchy-Schwarz inequality yields
  $\|x^2-y^2\|_{p/2}\le \|x-y\|_p\cdot \|x+y\|_p$,  for $x$ and $y\in L_p$, we may suppose that
$(x_n)$ is a  block basis of the Haar basis. Since the
sequence $x_n^2$ is positive and weakly convergent,
$$
\|x_n^2\|_1=\E[x_n^2]\to \E[v]= \|v\|_1=1. $$
  We can thus assume
that $\|x_n\|_2=1$ for all $n$.  We may choose a natural number $K$
such that $\|\E[v\mid\mathcal H_K]\|_{p/2}> M'^2 $ and 
by discarding the first few elements
of   $(x_n)$ we have that
 $x_n\in [h_k]_{k>K}$, for
all $n$. The $x_n$ are martingale differences with respect to a
subsequence $\F_n = \mathcal H_{k(n)}$ of the Haar filtration (with
$k(0)=K$).
Taking a further subsequence, we may suppose that
\begin{equation}
 \big\|\E[v\mid\mathcal F_{n-1}]^{1/2}-\E[x_n^2\mid
  \mathcal F_{n-1}]^{1/2}\big\|_\infty <2^{-n}, \text{ for all $n$.}\label{vxI}
\end{equation} 
 Because $(x_n)$ is a martingale 
  difference sequence, 
  we can apply Theorem~\ref{B-R-I} to conclude that
  \begin{equation*}
\Big\|\sum c_nx_n\Big\|_p \ge C_p^{-1} \Big\|\Big(\sum c_n^2
\E[x_n^2|\mathcal F_{n-1}]\Big)^{1/2}\Big\|_{p}\!=
\!C_p^{-1}\big\|  \|(c_n \E^{1/2}[x^2_n|\mathcal F_{n-1}]: n\!\in\!\nat)\|_{\ell_2}   \big\|_p.
\end{equation*}
If we use \eqref{vxI} and apply the triangle inequality in
$L_p(\ell_2)$ we obtain
\begin{align*}
\Big\|\sum c_nx_n\Big\|_p
&\ge C_p^{-1}\big\|  \|(c_n \E^{1/2}[x^2_n|\mathcal F_{n-1}]: n\in\nat)\|_{\ell_2}   \big\|_p \\
&\ge  C_p^{-1}\Big(\big\|  \|(c_n \E^{1/2}[v |\mathcal F_{n-1}]: n\in\nat)\|_{\ell_2}   \big\|_p 
  -\|(c_n 2^{-n}:n\in\nat)\|_{\ell_2}\Big)\\
&= C_p^{-1}\Big(\Big\|\big(\sum c_n^2
 \E[v|{\mathcal F}_{n-1}]\big)^{1/2}\Big\|_{p}\!-\!  \big(\sum c_n^2 2^{-2n}\big)^{1/2}\Big)
\ge \frac{M'-1}{C_p}\Big(\sum c_n^2\Big)^{1/2}.
\end{align*}

On the other hand, in $L_2$, the $x_n$ are orthogonal, whence
$$
\Big\|\sum c_nx_n\Big\|_2 =\Big(\sum c_n^2\Big)^{1/2}.
$$
Provided $M'$ is chosen large enough, we have $\|y\|_p\ge M\|y\|_2$
for all $y\in [x_n]$ as required.
\end{proof}

The next step is to show that we can choose our
$\ell_2$-subspaces to have $p$-uniformly integrable unit balls.
Recall that a subset $A$ of $L_p$ is said to be {\em $p$-uniformly
integrable} if, for every $\ep >0$ there exists $K>0$ such that
$\|x1_{[|x|>K]}\|_p<\ep$ for all $x\in A$. We shall need the
following standard martingale lemma.

\begin{lem}\label{pUI}
 Let $(x_n)$ be a martingale difference sequence that is $p$-uniformly integrable.  Then the
 set of linear combinations of the $x_n$'s with $\ell_2$-normalized coefficients is also
 $p$-uniformly
 integrable.
 \end{lem}
 \begin{proof}
 We assume that $\|x_n\|_2\le1$ for all $n$ and consider a vector $y$
 of the form $\sum_n c_nx_n$ with $\sum_n c_n^2= 1$, noting that
 $\|y\|_2^2 = \sum c_n^2 \|x_n\|_2^2\le 1$. Given $\ep>0$, we choose $K>\ep^{-1}$ such that
 $\|x_j 1_E\|_2<\ep$ for all $j$ whenever $\mathbb P(E)<K^{-1}$. 
 We consider the martingale $(y_n)$
 where $y_n = \sum_{j\le n}c_jx_j$ (thus $y= y_\infty$) and
 introduce the stopping time
 $$
 \tau = \inf\{n\in \mathbb N:|y_n|>K\}.
 $$
 By Doob's inequality $\mathbb P[\tau<\infty]\le K^{-1}\|y\|_1\le
 K^{-1}$.
 We note that if $\tau<\infty$, then $|y_\tau|\le K + |c_\tau x_\tau|$ so that
 $$
 |y| \le K + |y-y_\tau| + |c_\tau x_\tau1_{[\tau<\infty]}|.
 $$
 We shall estimate the $L_p$-norms of the second two terms.
 For the first of these, we note that $(y_k-y_{k\wedge \tau})$ is a
 martingale, so that ($ C$ only depends on $p$)

\begin{align*}
  \|y-y_\tau\|_p &\le C \Big\|\sum_n c_n^2 x_n^2
 1_{[\tau<n]} \Big\|_{p/2}^{1/2}\quad \text{[by\ the\ square\ function\ inequality]}\\
& \le C \Big(\sum c_n^2 \|x_n^2 1_{[\tau<n]}\|_{p/2}\Big)^{1/2}\quad\text{[by\
the\ triangle\ inequality\ in\ $L_{p/2}$]}\\
 &\le C \sup_n \|x_n1_{[\tau<\infty]}\|_p\qquad\big[\text{since\ $\sum c_n^2\le
 1$}\big]\\
 &\le C \ep\qquad\quad \text{[because\ $\mathbb P[\tau<\infty]\le
 K^{-1}$]}.\\
 \intertext{
For the second term we use the fact that the sets $[\tau=n]$ are
disjoint, so that}
\|c_\tau x_\tau1_{[\tau<\infty]}\|_p &= \Big\|\sum_nc_nx_n1_{[\tau=n]}\Big\|_p
  =\Big(\sum_n |c_n|^p \|x_n1_{[\tau=n]}\|_p^p\Big)^{1/p}
\le \sup_n \|x_n1_{[\tau<\infty]}\|_p \le \ep
\end{align*}
as before.  Thus,
$$\| y1_{[|y|>2K]}\|_p\le K{\mathbb P}^{1/p}\big[|y-y_\tau|\!+ \!|c_\tau x_\tau1_{[\tau<\infty]}|>K\big]+ (C+1)\ep
\le  2(1+C)\ep, $$
 which implies our claim.
 \end{proof}

\begin{lem}\label{L:5.4}
Let $Y$ be a subspace of $L_p$ ($p>2$), which is isomorphic to
$\ell_2$.  There is an infinite dimensional subspace $Z$ of $Y$ such
that the unit ball $B_Z$ is $p$-uniformly integrable.
\end{lem}
\begin{proof}
Let $(y_n)$ be a normalized sequence in $Y$ equivalent to the unit
vector basis of $\ell_2$. By the Subsequence Splitting Lemma (see,
for instance Theorem IV.2.8 of \cite{G-D}), we can write
$y_n=x_n+z_n$, where the sequence $(x_n)$ is $p$-uniformly
integrable, and the $z_n$ are disjointly supported. So $(x_n)$ and $(z_n)$ are weakly null. Taking a
subsequence, we may suppose that the $(x_n)$  is a  martingale
difference sequence, so that the set of all $\ell_2$-normalized linear
combinations $\sum c_nx_n$ is also $p$-uniformly integrable.

We now consider $\ell_2$-normalized blocks of the form
\begin{align*}
y_k'&= (N_k-N_{k-1})^{-1/2}\sum_{N_{k-1}<n\le N_k} y_n
= x_k'+z_k',
 \end{align*}
 where,
 $$ x'_k=   (N_k-N_{k-1})^{-1/2}\sum_{N_{k-1}<n\le N_k}
x_n   \text{ and }  z'_k=(N_k-N_{k-1})^{-1/2}\sum_{N_{k-1}<n\le N_k}
z_n.$$
 
Because the $z_n$ are disjointly supported in $L_p$ we have
$\|z'_k\|_p \le (N_k-N_{k-1})^{1/p-1/2}$, so we can choose the $N_k$
 such  that $\|z'_k\|_p<2^{-k}$. The sequence  $(x'_k)$, being $\ell_2$
normalized linear combinations of the $x_n$, are $p$-uniformly
integrable.  Hence the $y'_k$, which are small perturbations of the
$x'_k$, are also $p$-uniformly integrable.  Another application of
Lemma~\ref{pUI} yields the result.
\end{proof}

We are now ready for the proof of Theorem~\ref{SkinnyTed}.

\begin{proof}[Proof of  Theorem \ref{SkinnyTed}]
By  Lemmas   \ref{Thinell2} and  \ref{L:5.4} there exists, for each $M>0$, a subspace $Z_M$ of $X$,
isomorphic to $\ell_2$ with $p$-uniformly integrable unit ball,
such that
$$
\|y\|_p\ge M\|y\|_2
$$
for all $y\in Z_M$.  For a specified $\ep>0$, we shall choose
inductively $M_1<M_2<\cdots$ and define $Y_n= Z_{M_n}$, 
such that
\begin{equation}\label{wedge}
\||y_m|\wedge |y_n|\|_p\le \ep/n 2^{n},
\end{equation}
whenever $y_m\in B_{Y_m}$, $y_n\in B_{Y_n}$ and $m<n$.

To achieve this, we start by taking an arbitrary value for $M_1$,
say $M_=1$.  Recursively, if $M_1,\dots,M_n$ have been chosen,
we use the $p$-uniform integrability of $\bigcup_{m\le n} B_{Y_m}$
to find $K_n$ such that $\big\||y|-|y|\wedge K_n\big\|_p<\ep/(n+1)2^{n+2}$
whenever $y\in B_{Y_m}$ and $m\le n$. We now choose $M_{n+1}$ such
 that $M_{n+1}^2>K_n^{p-2}(n+1)^p2^{p(n+2)}\ep^{-p}$.

We need to check that \eqref{wedge} is satisfied, so let $y_{n+1}\in
B_{Y_{n+1}}$ and let $y_m\in B_{Y_m}$ with $m\le n$.  We have
\begin{equation*}
|y_m|\wedge |y_{n+1}| \le K_{n}\wedge |y_{n+1}| +(|y_m|-|y_m|\wedge
K_n)
\end{equation*}
and have chosen $K_n$ in such a way as to ensure that
$$
\big\||y_m|-|y_m|\wedge K_n \big\|_p<\ep/(n+1)2^{n+2}.
$$
For the first term, we note that
\begin{equation*}
\E[(K_n\wedge |y_{n+1}|)^p] \le \E[K_n^{p-2}|y_{n+1}|^2]
 = K_n^{p-2}\|y_{n+1}\|_2^2 \le K_n^{p-2}M_{n+1}^{-2},
\end{equation*}
which is smaller than $\ep^p(n+1)^{-p}2^{-p(n+2)}$, by our choice of $M_{n+1}$.

Now let $y_n\in S_{Y_n}$ for all $n\in\nat$.  We shall show that  the $y_n$'s are 
small perturbations of elements that are disjoint in $L_p$.  Indeed,
let us set
$$
y'_n = \text{sign}\,(y_n)\big(|y_n|-|y_n|\wedge\bigvee_{m\ne n} |y_m|\big).
$$
Then the  $y'_n$ are  disjointly supported and from \eqref{wedge}
\begin{equation*}
\|y_n-y'_n\|_p = \Big\||y_n|\wedge\bigvee_{m\ne n} |y_m|\Big\|_p
  \le \sum_{m\ne n}\big \||y_n|\wedge |y_m|\big\|_p
  \le (n-1)\ep/n2^{n} + \sum_{m>n} \ep/m2^m<\ep/2^n.
  \end{equation*}
Standard manipulation of inequalities now shows us that the closure
of the sum $\sum_n Y_n$ in $L_p$ is almost an $\ell_p$-sum.  Indeed,
\begin{align*}
(1-2\ep) \Big(\sum |c_n|^p\Big)^{1/p} &\le \Big(\sum
|c_n|^p\|y'_n\|_p^p\Big)^{1/p}-\ep\Big(\sum |c_n|^p\Big)^{1/p}\\
&= \Big\|\sum c_ny'_n\Big\|_p -\ep\Big(\sum |c_n|^p\Big)^{1/p}\\
&\le \Big\|\sum c_n y_n\Big\|_p\\
&\le \Big\|\sum c_n y'_n\Big\|_p + \ep\Big(\sum |c_n|^p\Big)^{1/p}
\le (1+\ep)\Big(\sum |c_n|^p\Big)^{1/p}.
\end{align*}
At this point in the proof, we have obtained subspaces $Y_n$ of $X$,
each isomorphic to $\ell_2$ such that the closed linear span
$\overline{\sum_n Y_n}$ is almost isometric to $(\bigoplus Y_n)_p$.
By stability (\cite{KM} or \cite{AO}) we can take, for each $n$, a
subspace $X_n$ of $Y_n$ which is $(1+\ep)$-isomorphic to $\ell_2$.
In this way we obtain a subspace of $X$ which is almost isometric to
$\ell_p(\ell_2)$. \end{proof}

The last part of the claim of Theorem B, namely that  we can pass to  a further subspace
 of $X$ which is still $(1+\theta)$-isomorphic to $\ell_p(\ell_2)$ and, moreover,
  complemented in $L_p$ follows from our results in the next section.
  G.~Schechtman \cite{S2} showed us that if one is not concerned with minimizing the 
norm of the projection, then there is a short argument that gives a complemented
copy of $\ell_p(\ell_2)$. We thank him for allowing us to present it here.

\begin{prop}  Let  $X\subset L_p$  be isomorphically equivalent
  to $\ell_p(\ell_2)$.
Then there is a subspace $Y$ of $X$ which is isomorphic to $\ell_p(\ell_2)$  and complemented in $L_p$.
\end{prop}

\begin{proof}
Let $\{x(m,n):m,n\in\nat\}\subset X$ be a normalized  basis of  $X$ equivalent to the usual
 unconditional
basis of $\ell_p(\ell_2)$, i.e. there is a constant $C\ge 1$ so that
$$\Big\|\sum_{m,n\in\nat} a(m,n) x(m,n)\Big\| \Csim \Big( \sum_{m\in\nat}\Big(\sum_{n\in\nat} a(m,n)^2\Big)^{p/2}\Big)^{1/p}
\text{ for all }(a(m,n))\in c_{00}(\nat^2).$$

 In  \cite{PR} it was shown that  for any $C>1$  there  is a $g_p(C)<\infty$ so that
 every subspace $E$ of $L_p$, which is $C$ isomorphic to $\ell_2$, is $g_p(C)$ complemented
  in $L_p$.
For $m\in\nat$ let $P_m: L_p\to [(x(m,n):n\in\nat]$ be a projection 
 of norm at most $g_p(C)$.  We can write 
 $$P_m(x)=\sum_{n\in\nat} x^*(m,n)(x) x(m,n)\text{ for } x\in L_p,$$
where $(x^*{(m,n)}:n\in\nat)$ is a weakly null sequence  in $L_q$, $\frac1p+\frac1q=1$,
and biorthogonal to $x(m,n):n\in\nat)$.
By passing to subsequences, using a diagonal argument, and perturbing we may assume that
there is a blocking  $(H(m,n):m,n\in\nat)$ of the Haar basis of $L_p$,
 in some order, so that $x(m,n)\in H(m,n)$ and $x^*(m,n)\in H^*(m,n)$, for $m,n\in\nat$, where
  $(H^*(m,n))$ denotes the blocking  of the Haar basis in $L_q$
  which corresponds to $(H(m,n))$
 
 We will show that the operator
 $$P: L_p\to L_p,\quad x\mapsto \sum_{m,n\in\nat} x^*(m,n)(x) x(m,n),$$
 is bounded and, thus, it is a bounded projection onto $[x(m,n):m,n\in\nat]$.
 
For $y=\sum_{m,n\in\nat} y(m,n)$, with $y(m,n)\in H(m,n)$, if $m,n\in\nat$, we deduce that
\begin{align*}
\|P(y)\|&=\Big\|\sum_{m\in\nat}\sum_{n\in\nat} x^*(m,n)(y(m,n)) x(m,n)\Big\|\\
         &\le  C\Big(\sum_{m\in\nat}\Big(\sum_{n\in\nat} (x^*(m,n)(y(m,n)))^2\Big)^{p/2}\Big)^{1/p}\\
           &\le C^2\Big(\sum_{m\in\nat} \|P_m(y_m)\|^p\Big)^{1/p}\le  C^2 g_p(C)
           \big(\sum_{m\in\nat} \|y_m\|^p\Big)^{1/p}
\end{align*}where $y_m=\sum_{n\in\nat} y(m,n)$  for $m\in\nat$.

The Haar basis is unconditional in $L_p$, and  if we denote 
the unconditional constant  in $L_p$ by $U_p$ we deduce from
 Proposition \ref{prop2.1} that
$\|y\|\ge U_p ^{-1}      (\sum_{m\in\nat}\|y_m\|^p)^{1/p}$,
which implies our claim.
\end{proof}

\begin{remark} G.~Schechtman \cite{S2} has also proved,
by a more complicated argument, that if $X\subset L_p$, $1<p<2$, is an
isomorph of $\ell_p(\ell_2)$ then $X$  contains a copy of $\ell_p(\ell_2)$ which is complemented
in $L_p$.
\end{remark}
Let us now deduce the statement of Corollary D.
\begin{proof}[Proof of Corollary D]
First assume that $X$ embeds into $\ell_p\oplus\ell_2$. Note that every
weakly null sequence $(x_n)$ can be turned into a weakly null tree $(x_\alpha)$, whose 
branches are exactly the subsequences of $(x_n)$ (put
 $x_{(n_1,n_2,\ldots n_\ell)}=x_{n_\ell}$ for 
  $(n_1,n_2,\ldots n_\ell)\in T_\infty$). This fact, together with the remarks at the beginning 
 of the proof of Theorem A (about the existence of $K$), show that  condition (b)  of Theorem A for a subspace $X$ of $L_p$  implies, that
   there exists a $K\ge 1$,  so that
 every weakly null sequence  in $S_X$ admits a subsequence  $(x_i)$ satisfying
   for all scalars $(a_i)$ condition  \eqref{eqA1} in (b) of Theorem A.
   
   Conversely, assume that $X$ does not embed into $\ell_p\oplus\ell_2$. Then Propositions 
 \ref{K-P}  and \ref{P:4.2} together with Theorem A imply that condition (B)  of 
  Proposition \ref{K-P} is satisfied.  Using now  Lemma \ref{Thinell2}, we can find for every
   $M<\infty$ a subspace $Y$ of $X$ which is isomorphic to $\ell_2$, so that
   $\|\cdot\|_p\ge M\|\cdot\|_2$ on $Y$. This implies that there cannot be a $K\ge 1$,  so that
 every weakly null sequence  in $S_X$ admits a subsequence  $(x_i)$ satisfying
   \eqref{E:1.4}.
\end{proof}

\section{Improving the embedding via random measures}\label{S:6}

We shall give a quick review of what we need from the theory of
stable spaces and random measures.  We shall then obtain the
optimally complemented embeddings of $\ell_p(\ell_2)$.

We start this section by recalling some facts about random measures
and their relation to types on $L_p$. The introductory part is valid for  $1<p<\infty$. Later we will restrict 
ourselves again to the case $p>2$. As far as possible, we shall
follow the notation and terminology of \cite{A}; for the theory of
types and stability we refer the reader to \cite{KM} (or \cite{AO}).
The lecture notes of Garling \cite{Garling} is one of the few works
where the connection between random measures and types on function
spaces is explicitly considered.

We shall denote by $\mathcal P$ the set of probability measures on
$\real$ which is a Polish space for its usual topology.  This
topology, often called the ``narrow topology'', can be thought of as
the topology induced by the weak* topology $\sigma(\mathcal
C_b(\real)^*, \mathcal C_b(\real))$.  

A {\it random measure} on $(\Omega, \Sigma,\mathbb P)$ is a
mapping $\xi:\omega\mapsto \xi_\omega;\Omega\to \mathcal P$ which is
measurable from $\Sigma$ to the Borel $\sigma$-algebra of
$\mathcal P$.  The set of all such random measures is denoted by
$\mathcal M$ and is a Polish space when equipped with what Aldous
calls the {\it wm-topology}. Sequential convergence for this
topology can be characterized by saying that $\xi^{(n)}\wmto \xi$ if
and only if
$$
\E\left[1_F\int_\real f(t)\text d\xi^{(n)}(t)\right] \to
\E\left[1_F\int_\real f(t)\text d\xi(t)\right],
$$
for all $F\in \Sigma$ and all $f\in \mathcal C^{\text b}(\real)$.  In
interpreting the expectation operator in the above formula (and in
similar expressions involving ``implicit'' $\omega$'s) the reader
should bear in mind that $\xi$ is random. If we translate the
expectation into integral notation,
$$
\E\left[1_F\int_\real f(t)\text d\xi(t)\right]
\text{ becomes }
 \int_F \int_\real f(t)\, \text d\xi_\omega(t)\,\text d\mathbb P(\omega).
$$
It is sometimes useful to use the notation $\xi_F$, when $F$ is a
non-null set in $\Sigma$ for the probability measure given by
$$
\int_\real f(t)\, \text d\xi_F(t) = \mathbb
P(F)^{-1}\E[1_F\int_\real f(t)\, \text d\xi(t)]\quad(f\in \mathcal
C_0(\real)).
$$

The usual convolution operation on $\mathcal P$ may be extended to
an operation on $\mathcal M$ by defining $\xi*\eta$ to be the
random measure with $(\xi*\eta)_\omega= \xi_\omega*\eta_\omega$.
Garling (Proposition 8 of \cite{Garling}) observes that this
operation is separately continuous for the wm topology. This result
is also implicit in Lemma 3.14 of \cite{A}. We may also introduce a
``scalar multiplication'': when $\xi\in\mathcal M$ and $\alpha$ is
a random variable, we define the random measure $\alpha. \xi$ by
setting
$$
\int f(t) \,\text d(\alpha.\xi)(t) = \int_\real f(\alpha t)\,\text
d\xi(t) \quad(f\in \mathcal C^{\text b}(\real)).
$$

Every random variable $x$ on $(\Omega,\Sigma,\mathbb P)$ defines a
random (Dirac) measure $\omega\mapsto \delta_{x(\omega)}$. Aldous \cite[after Lemma 2.14]{A}
has remarked that (provided the probability space
$(\Omega,\Sigma,\mathbb P)$ is atomless)  these $\delta_x$ form a
wm-dense subset of $\mathcal M$. While we do not need this fact
here, it may be helpful to note that the definition given above of
$\alpha. \xi$ is so chosen that $\delta_{\alpha x_n} \wmto \alpha.
\xi$ whenever $\delta_{x_n}\wmto \xi$.  The $L_p$-norms extend to
wm-lower semicontinuous $[0,\infty]$-valued functions $|\cdot|_p$ on
$\mathcal M$, defined by
$$
|\xi|_p = \E\Big[\int_\real |t|^p \,\text d\xi(t)\Big]^{1/p}.
$$
We shall write $\mathcal{M}_p$ for the set of all $\xi$ for which
$|\xi|_p$ is finite.

As a special case of the characterization of wm-compactness by the
condition of ``tightness'' we note that a subset of $\mathcal M_p$
which is bounded for $|\cdot|_p$ is wm-relatively compact.  In
particular, if $(x_n)$ is a sequence that is bounded in $L_p$ then
there is a subsequence $(x_{n_k})$ such that $\delta_{x_{n_k}}\wmto
\xi$ for some $\xi\in\mathcal M_p$. 
 If $(x_n)$ is, moreover, $p$-uniformly integrable,  an easy truncation argument
 shows that 
$$\lim_{n\to\infty} \|x_n\|_p=
\lim_{n\to\infty}\E\Big(\int |t|^p d\delta_{x_n}(t)\Big)=
\E\Big(\int |t|^p d\xi(t)\Big).$$
For a subspace $X$ of $L_p$ we
write $\mathcal M_p(X)$ for the set of all $\xi$ that arise as
wm-limits of sequences $(\delta_{x_n})$ with $(x_n)$ an
$L_p$-bounded sequence in $X$. It is an easy consequence of separate
continuity that $\mathcal M_p(X)$ is closed under the convolution
operation $*$ (c.f. the proof of \cite[Proposition 3.9]{A}).

We recall that a function $\tau:X\to \real$ on a (separable) Banach
space $X$ is called a {\em type} if there is a sequence $(x_n)$ in
$X$ such that, for all $y\in X$,
$$
\|x_n+y\|\to \tau(y)\quad\text{as\ }n\to \infty.
$$
The set of all types on $X$ is denoted $\mathcal T_X$ and is a
locally compact Polish space for the {\em weak } topology; this
topology may be characterized by saying that $\tau_n\wto \tau$ if
$\tau_n(y)\to \tau(y)$ for all $y\in X$.  If we introduce, for each
$x\in X$, the {\em degenerate type} $\tau_x$ defined by
$$
\tau_x(y)= \|x+y\|,
$$
then $\mathcal T_X$ is the w-closure of the set of all $\tau_x$. We
introduce a ``scalar multiplication'' of types, defining
$\alpha.\tau$, for $\alpha\in \mathbb R$ and $\tau\in \mathcal T_X$
by setting
$$
\alpha. \tau = \text{w-lim}\, \tau_{\alpha x_n}\quad\text{when}
\quad \tau= \text{w-lim}\,\tau_{x_n}.
$$

A Banach space $X$ is {\em stable} if, for $x_m$ and $y_n$ in $X$,
we have
$$
\lim_{m\to \infty}\lim_{n\to \infty} \|x_m+y_n\|=
\lim_{n\to \infty}\lim_{n\to \infty} \|x_m+y_n\|,
$$
whenever the relevant limits exist.  All $L_p$-spaces ($1\le
p<\infty$) are stable \cite{KM}.

Stability of a Banach space $X$ permits the introduction of a
(commutative) binary operation $*$ on $\mathcal T_X$, defined by
$$
\tau*\upsilon(z) = \lim_{m\to \infty} \lim_{n\to
\infty}\|x_m+y_n+z\|
$$
when $\tau=\text{w-lim}\,\tau_{x_m}$ and $\upsilon =
\text{w-lim}\,\tau_{y_n}$.

A type $\tau\in \mathcal T_X$ is said to be an $\ell_q$-type if
$$
(\alpha.\tau)*(\beta.\tau) = (|\alpha|^q+|\beta|^q)^{1/q}.\tau
$$
for all real $\alpha, \beta$.  The big theorem of \cite{KM} shows
first that on every stable space there are $\ell_q$-types for some
value(s) of $q$, and secondly that the existence of an $\ell_q$ type
implies that the space has subspaces almost isometric to $\ell_q$.
In fact the proof of Th\'eor\`eme III.1 in \cite{KM} proves
something slightly more than the existence of such a subspace.  We
now record the statement we shall need.

\begin{prop}\label{StableSubspace}
Let $X$ be a stable Banach space, let $1\le q<\infty$ and let
$(x_n)$ be a sequence in $X$ such that $\tau_{x_n}$ converges to an
$\ell_q$-type $\tau$ on $X$. Then there is a subsequence $(x_{n_k})$
such that $\tau_{z_n}$ converges to $\tau$ for every
$\ell_q$-normalized block subsequence $(z_n)$ of $(x_{n_k})$.
\end{prop}

The results of \cite{KM} extended, and gave an alternative approach
to  the theorem of \cite{A}, which obtained $\ell_q$'s in subspaces
of $L_1$ using random measures.  We shall need elements from both
approaches.  The link is provided by the following lemma, for which
we refer the reader to the final paragraphs of \cite{Garling}. We
shall write $\mathcal T_p$ for $\mathcal T_{L_p}$ and, when $X$ is
a subspace of $L_p$, we shall write $\mathcal T_p(X)$ for the weak
closure in $\mathcal T_p$ of the set of all $\tau_x$ with $x\in X$.

\begin{lem}
Let $(x_n)$ be a bounded sequence in $L_p$ and suppose that
$\delta_{x_n}\wmto \xi$ in $\mathcal M$.  Suppose further that
$\|x_n\|_p\to \alpha$ as $n\to \infty$.  Then, for all $y\in L_p$
$$
\|x_n+y\|^p_p \to \E\left[\int_\real |y+t|^p\, \text d\xi(t)\right ]
+ \beta^p,
$$
where the non-negative constant $\beta$ is given by
$$
\alpha^p = \|\xi\|_p^p + \beta^p.
$$
The sequence $(x_n)$ is $p$-uniformly integrable if and only if
$\beta=0$.
\end{lem}

We thus have the following formula showing how the type $\tau=\lim
\tau_{x_n}\in \mathcal T_p$ is related to the random measure $\xi =
\text{wm-}\lim \delta_{x_n}\in \mathcal M_p$ and the index of
$p$-uniform integrability $\beta$.
\begin{equation}\label{TypeRep}
\tau(y)^p = \E\left[\int_\real |y+t|^p\text d\xi(t)\right] +
\beta^p.
\end{equation}
If $q<p$ then a sequence $(x_n)$ as above in $L_p$ can be thought of
as a sequence in $L_q$.  If we wish to distinguish the type
determined  on $L_q$ from the type on $L_p$, we use
superscripts.  Of course,
$$
\tau^{(q)}(y)^q = \E\left[\int_\real |y+t|^q\text d\xi(t)\right],
$$
with no ``$\beta$'' term, because an $L_p$-bounded sequence is
$q$-uniformly integrable.

The * operations on $\mathcal T_p$ and on $\mathcal M_p$ are
related by the following lemma, also to be found in \cite{Garling}.

\begin{lem}\label{StarRep}
Let $\tau_1$ and $\tau_2$ be types on $L_p$ represented as
\begin{equation*}
\tau_1(y)^p = \E\left[\int_\real |y+t|^p\text d\xi_1(t)\right] +
\beta_1^p\text{ and }
\tau_2(y)^p = \E\left[\int_\real |y+t|^p\text d\xi_2(t)\right] +
\beta_2^p.
\end{equation*}
Then
$$
(\tau_1*\tau_2)(y)^p =  \E\left[\int_\real |y+t|^p\text
d(\xi_1*\xi_2)(t)\right] +\beta_1^p+  \beta_2^p.
$$
\end{lem}

It has been noted already in the literature (e.g. \cite{Garling})
that the representation given in \eqref{TypeRep} is not in
general unique. However, for most values of $p$, it is, as we now
show.

\begin{prop}\label{UniqueRM}
Let $1\le p<\infty$ and assume that $p$ is not an even integer.  In
the representation of a type $\tau$ on $L_p$ by the formula
\eqref{TypeRep} the random measure $\xi$ and the constant $\beta$ are
uniquely determined by $\tau$. If $(x_n)$ is any sequence in $L_p$
with $\tau_{x_n}\wto \tau$ we have $\delta_{x_n}\wmto \xi$ and
$\inf_M\lim_{n\to \infty}\|x_n1_{[|x_n|\ge M]}\|_p= \beta$.
\end{prop}
\begin{proof} Suppose that $\xi,\beta$ and $\xi',\beta'$ yield the
same type $\tau$.  For any non-null $E\in \Sigma$ and any real number
$u$, we consider $\tau(y)$ where $y= u1_E\in L_p$ to obtain
$$
\E\left[\int_\real |t+u1_E|^p \text d\xi(t)\right] +\beta^p =
\E\left[\int_\real |t+u1_E|^p \text d\xi'(t)\right] +\beta'^p,
$$
or, equivalently,
$$
\int_\real |t+u|^p \text d\xi_E(t) = \int_\real |t+u|^p \text
d\xi'_E(t) + \alpha^p,
$$
where
$$\Pr(E)\alpha^p = \beta'^p - \beta^p + \E\left[1_{\Omega\setminus E}\int
|t|^p\text d\xi'(t)-1_{\Omega\setminus E}\int |t|^p\text
d\xi(t)\right].
$$
By the Equimeasurability Theorem (cf.\cite[page 903]{KK}), $\alpha=0$ and the
measures $\xi_E$ and $\xi'_E$ are equal.  Since this is true for all
$E$, $\xi=\xi'$.

Now let $(x_n)$ be any sequence with $\tau_{x_n}\wto \tau$. By the
uniqueness that we have just proved, the only cluster point of the
sequence $\delta_{x_n}$ in $\mathcal M$ is $\xi$.  Since (by
$L_1$-boundedness) $\{\delta_{x_n}: n\in \nat\}$ is relatively
wm-compact in $\mathcal M$, it must be that $\delta_{x_n}\wmto \xi$.
\end{proof}

We have already noted that $\mathcal M_p(X)$ is closed under $*$
when $X$ is a subspace of $L_p$.  The next proposition, which is
closely related to  that of \cite[Proposition 3.9]{A}, shows that under
appropriate conditions $\mathcal M_p(X)$ is wm-closed.

\begin{prop}\label{MpClosed} Let $1\le p<\infty$ and let $X$ be a
subspace of $L_p$ with no subspace isomorphic to $\ell_p$. Then
$\mathcal M_p(X)$ is wm-closed in $\mathcal M$.
\end{prop}
\begin{proof} The hypothesis implies that the $L_p$-norm is
equivalent to the $L_1$-norm on $X$, so that we may regard $X$ as a
(reflexive) subspace of $L_1$. Aldous  \cite[Lemma 3.12]{A} shows (by a
straightforward uniform integrability argument) that
$\xi\mapsto|\xi|_1$ is wm-continuous and finite on $\mathcal D$,
where $\mathcal D$ is the wm-closure of $\{\delta_x:x\in X\}$. Thus
every $\xi$ in $\mathcal D$ is in the wm-closure of an
$L_1$-bounded subset of $X$, and hence, by equivalence of norms, in
$\mathcal M_p(X)$. \end{proof}

To finish this round-up of types and random measures, we need to
mention the connection between $\ell_2$-types and the normal
distribution  (a special case of the connection between
$\ell_q$-types and symmetric stable laws).  We write $\gamma$ for
the probability measure (or {\em law}) of a standard $\mathcal
N(0,1)$ random variable.  If $\sigma$ is a non-negative random
variable then $\sigma.\gamma$ is a random measure (a normal
distribution with random variance).  Provided $\sigma\in L_p$ this
random measure defines a type on $L_p$ by
$$
\tau(y)^p = \E\left [\int_\real |y+t|^p \,\text
d(\sigma.\gamma)(t)\right] = \E\left[\int_\real |y+\sigma t|^p\,
\text d\gamma(t)\right].
$$
Now it is a property of the normal distribution that
$(\alpha.\gamma)*(\beta.\gamma)= (\alpha^2+\beta^2)^{1/2}.\gamma$
for real $\alpha,\beta$. By Lemma~\ref{StarRep}, this allows us to see
that $\tau$ is an $\ell_2$-type on $L_p$.

We are finally ready to return to the main subject matter of this
paper.

\begin{lem}\label{VXlem}
Let $X$ be a subspace of $L_p$, with $p>2$, and let $v$ be a non-zero
element of $ L_{p/2}$. The following are equivalent:
\begin{enumerate}
\item $v\in V(X)$;
\item there exists $\xi\in \mathcal M_p(X)$ such that $\int_\real t\,
\text d\xi=0$ and $\int_\real t^2\,\text d\xi=v$ almost surely;
 \item $\sqrt v.\gamma\in \mathcal M_p(X)$.
\end{enumerate}
\end{lem}
\begin{proof}
 We start by assuming (1). Let $(x_n)$ be a weakly null sequence
in $X$ such that $(x_n^2)$ converges weakly to $v$ in $L_{p/2}$.
Replacing $(x_n)$ with a subsequence, we may suppose that
$\delta_{x_n}\wmto \xi$ for some $\xi\in \mathcal M_p(X)$. Since the
sequence $(x_n)$ is $L_p$-bounded, it is 2-uniformly integrable and
so
\begin{align}
\E\left[1_E\int_\real t \text
d\xi(t)\right]&=\lim\E\left[1_Ex_n\right] =0\quad\text{and}\\
\E\left[1_E\int_\real t^2 \text d\xi(t)\right]&=\lim
\E\left[1_Ex_n^2\right] = \E\left[1_Ev\right],
\end{align}
for all $E\in \Sigma$. This yields (2).

We now assume (2). Let  $(x_n)$ be an $L_p$-bounded sequence in $X$
such that $\delta_{x_n}$ is wm-convergent to  $ \xi$.
Since $\int_{\mathbb R} d\xi(t)=0$ a.s. it follows that $(x_n)$ is weakly null and
 since $\xi \ne
\delta_0$, $\|x_n\|_2$ does not tend to zero.  
 By \cite{KP}, it follows that $X_0$, the closed linear span of  a subsequence of $(x_i)$, is isomorphic to 
  $\ell_2$.
The assumption about $\xi$ is that, for almost all
$\omega$, the probability measure $\xi_\omega$ is the law of a
random variable with mean 0 and variance $v(\omega)$.

By the Central Limit Theorem
$$
n^{-1/2}.\underbrace{(\xi_\omega*\xi_\omega*\cdots*\xi_\omega)}_{\text{$n$\
terms}}
$$
tends to $\sqrt {v(\omega)}.\gamma$ for all such $\omega$.  So in
$\mathcal M$ we have
$$
n^{-1/2}.(\xi*\xi*\cdots*\xi)\wmto \sqrt v.\gamma.
$$
Since $\mathcal{M}_p(X_0)$ is closed under convolution and is closed
in the wm-topology (by Proposition~\ref{MpClosed}), we see that
$\sqrt v.\gamma\in \mathcal{M}_p(X_0)\subseteq\mathcal{M}_p(X)$.

Finally, if we assume (3) we may take $(x_n)$ to be an $L_p$-bounded
sequence in $X$ such that $\delta_{x_n}\wmto \sqrt v.\gamma$.
Calculations like those used in the proof of (1) $\implies$ (2),
justified by 2-uniform integrability, show that $(x_n)$ is weakly
null and that $x_n^2$ tends weakly to $v$.
\end{proof}

We shall say that a sequence $(y_n)$ in $L_p$ is a {\it stabilized
$\ell_2$ sequence} with limiting conditional variance $v$ if, for
every $\ell_2$ normalized block subsequence $(z_n)$ of $(y_n)$, the
following are true:
\begin{align}
\delta_{z_n} &\wmto \sqrt v.\gamma \text{\ as\ }n\to \infty;\label{Sell2-1}\\
\|z_n\|_p &\to \gamma_p \|\sqrt v\|_p \text{\ as\
}n\to\infty.\label{Sell2-2}
\end{align}
(Recall that $\gamma_p = \|x\|_p$,  where $x$  is a symmetric $L_2$
normalized Gaussian random variable).
 For $p$ not an even integer, it is not hard to establish the
existence of such sequences using Proposition~\ref{StableSubspace} and
Proposition~\ref{UniqueRM}.  The proof of the next proposition avoids the
irritating problem posed by non-unique representations, by switching
briefly to the $L_1$-norm.

\begin{prop}\label{P:stabblock}
Let $X$ be a closed subspace of $L_p$ ($p\!>\!2$) and let $v$ be a non-zero
element of $V(X)$.  Then there exists a stabilized $\ell_2$ sequence in
$X$ with limiting conditional variance~$v$.
\end{prop}
\begin{proof}
By Lemma~\ref{VXlem} the random measure $\sqrt v.\gamma $ is in $\mathcal
M_p(X)$.  Let $(x_n)$ be a bounded sequence in $X$ with
$\delta_{x_n}\wmto \sqrt v.\gamma$. For the moment, think of the
$x_n$ as elements of $L_1$ and consider the types $\tau^{(1)}_{x_n}$
defined on $L_1$.  By $L_p$-boundedness, the sequence $(x_n)$ is
uniformly integrable, so the  sequence $(\tau^{(1)}_{x_n})$
converges weakly to the $\ell_2$-type $\tau^{(1)}$, where
$$
\tau^{(1)} (y) = \E\left[\int|y+\sqrt vt|\text d\gamma(t)\right].
$$
By Proposition~\ref{StableSubspace} we may replace $(x_n)$ by a subsequence in
such a way that $\tau^{(1)}_{z_n}\wto \tau^{(1)}$ for every
$\ell_2$-normalized block subsequence $(z_n)$.  By Proposition~\ref{UniqueRM} we
have $\delta_{z_n} \wmto \sqrt v.\gamma$ for all such $(z_n)$.

We now return to the $L_p$-norm, for which we can assume, after passing to
 a subsequence, if necessary, that $(x_n)$ is
equivalent to the unit vector basis of $\ell_2$.  By stability of
$L_p$ there is an $\ell_2$-normalized block subsequence $(y_n)$ such
that $\tau^{(p)}_{y_n}\wto \tau^{(p)}$ for some $\ell_2$-type
$\tau^{(p)}$ on $L_p$.  Moreover, by Proposition~\ref{StableSubspace} we can
arrange that $\tau^{(p)}_{z_n}\wto \tau^{(p)}$ for every further
such $\ell_2$-normalized block subsequence $(z_n)$. By \eqref{TypeRep} we have
$$
\tau^{(p)}(y)^p = \E\left[\int_\real |y+\sqrt v t|^p\text
d\gamma(t)\right] + \beta^p,
$$
for some non-negative constant $\beta$.  Now $\tau^{(p)}$ is an
$\ell_2$-type, so $\tau^{(p)}*\tau^{(p)}= \sqrt 2.\tau^{(p)}$. That
is to say
$$
(\tau^{(p)}*\tau^{(p)})(y)^p = \E\left[\int_\real |y+\sqrt{2 v}
t|^p\text d\gamma(t)\right] + (\sqrt2\beta)^p.
$$
On the other hand, by Lemma~\ref{StarRep},
\begin{equation*}
(\tau^{(p)}*\tau^{(p)})(y)^p = \E\left[\int_\real |y+\sqrt{v}
t|^p\text d(\gamma*\gamma)(t)\right] + 2\beta^p
= \E\left[\int_\real |y+\sqrt{2 v} t|^p\text d\gamma(t)\right] +
2\beta^p.
\end{equation*}
Since $p\ne 2$, we are forced to conclude that $\beta=0$.

To sum up, for every $\ell_2$-normalized block subsequence $(z_n)$
of $(y_n)$ we have, first of all, $\delta_{z_n}\!\!\wmto\!\!\sqrt v
\gamma$, since the $z_n$ are normalized blocks of $(x_n)$.  But also
$$
\|z_n\|_p \to \tau^{(p)}(0) = \E\left[ \int_\real |\sqrt vt|^p\text
d\gamma\right]^{1/p} = \gamma_p \|\sqrt v\|_p.
$$
\end{proof}

\begin{thm}\label{T:altB}
Let $X$ be a subspace of $L_p$ ($p>2$) and assume that (B) of 
Proposition~\ref{K-P} holds.  Then, for every $\theta>0$, there is a subspace
$Y$ of $X$ which is $(1+\theta)$-isomorphic to $\ell_p(\ell_2)$ and
a projection $P$ from $L_p$ onto $Y$ with $\|P\|\le (1+\theta)\gamma_p$.
\end{thm}
\begin{remark} The fact that Theorem \ref{T:altB} is the optimal result concerning the norm of a
projection onto a copy of $\ell_p(\ell_2)$ follows from \cite[Theorem 5.12]{GLR},
where it was show that $L_p$ contains subspaces isometric to $\ell_2$ which are
$\gamma_p$ complemented.
\end{remark}
\begin{proof}
Let $\ep\in (0,1)$ be fixed and, for $m\in\nat$, let $v_m\in V(X)$, together with
disjoint sets $A_m\in \Sigma$,  $A_m\subset \supp(v_m)$, be chosen so that
$\|v_m^{1/2}1_{A_m}\|_p = 1$ and 
$\|v_m^{1/2}\|^p_p < 1+\ep^p2^{-(m+2)p}$. Using Proposition \ref{P:stabblock} choose
for each $m$  a stabilized $\ell_2$-sequence
$(x_n^{(m)})_{n\in \nat}$ in $X$ with limiting conditional variance
$v_m$. By \eqref{Sell2-1} we have
\begin{equation*}
\liminf_{n\to\infty} \E[|y_n|^p1_{A_m}] \ge \gamma_p^p\quad\text{ and }
\liminf_{n\to\infty} \E[y_n^2v_m^{\frac p2-1}1_{A_m}] \ge 1 
\end{equation*}
and, by \eqref{Sell2-2},
$$ \lim_{n\to\infty} \E[|y_n|^p] = \gamma_p^p\|\sqrt v_m^{1/2}\|_p^p<
 \gamma_p^p(1+\ep^p2^{-(m+2)p}),
$$
 for all $\ell_2$-normalized block subsequences
$(y_n)$ of $(x^{(m)}_n)$. By relabeling the sequence $(x_n^{(m)})$,
starting at a suitably large value of $n$, we may suppose that the
following hold for all $\ell_2$-normalized linear combinations $y$
of the $x^{(m)}_n$:
\begin{align}
\|y1_{A_m}\|_p^p &\ge (1-\ep 2^{-(m+2)p})\gamma_p^p\label{Ym1}\\
\E\left[y^2v_m^{\frac p2-1}1_{A_m}\right] &\ge 1-\ep2^{-m-1}\label{Ym2}\\
\|y\|_p^p &\le (1+\ep 2^{-(m+2)p})\gamma_p^p\label{Ym3}.
\end{align}
Of course, \eqref{Ym1} and \eqref{Ym3} imply that the closed linear span
$Y_m=[x^{(m)}_n]_{n\in \nat}$ is almost isometric to $\ell_2$;
indeed, by homogeneity,  they yield
$$
(1-\ep 2^{-(m+2)p})^{1/p}\gamma_p(\sum c_n^2)^{1/2} \le \|y\|_p \le
(1+\ep 2^{-(m+2)p})^{1/p}\gamma_p(\sum c_n^2)^{1/2},
$$
when $y=\sum c_n x_n^{(m)}\in Y_m$.

Moreover, from the same inequalities we obtain
\begin{equation}
\|y-y1_{A_m}\|_p \le \ep 2^{-m}\|y\|_p\label{yy'}\text{ for all $y\in Y_m$.} 
\end{equation}
 If, for each $m\in\nat$, $y_m$ is an
element of $S_{Y_m}$ then $y_m'=y_m1_{A_m}$ are disjointly supported
and are small perturbations of the $y_m$. As in the proof of Theorem
\ref{SkinnyTed}, we see that, by an appropriate choice of $\ep$, we
can arrange for the closure of $\sum_m Y_m$ in $X$ to be
$(1+\theta)$-isomorphic to $\ell_p(\ell_2)$.  We are now ready to
show that the subspace $Y= \overline{\sum_m Y_m}$ is complemented in
$L_p$.  We shall do this by combining the disjoint perturbation
procedure used above with a standard ``change-of-density'' argument.

For each $m$ let $\phi_m= v_m^{p/2}1_{A_m}$; thus $\|\phi_m\|_1=1$.
Let $\Phi_m:L_p\to L_p(\phi_m)$ be defined by
$$
\Phi_m(f) = 1_{A_m} \phi_m^{-1/p}f,
$$
which is well defined since $A_m\subset\supp(v_m)$, and observe that
$$
\|\Phi_m(f)\|_{L_p(\phi_m)} = \|f1_{A_m}\|_p.
$$
Let $J_m: L_p(\phi_m)\to L_2(\phi_m)$ be the standard inclusion and
let $I_m:Y_m\to L_p$ be the natural embedding.  We note that for
$y\in Y_m$
\begin{equation*}
\|J_m\Phi_mI_my\|_{L_2(\phi_m)}^2 = \E[y^2 \phi_m^{-2/p}
\phi_m1_{A_m}]
=\E[y^2v_m^{\frac p2-1}1_{A_m}]
\ge (1-\ep2^{-m})^2\gamma_p^{-2}\|y\|_p^2,
\end{equation*}
by \eqref{Ym2}, \eqref{Ym3} and homogeneity. 
So if $W_m$ is the image
$$
W_m = J_m\Phi_mI_m[Y_m]
$$
then $W_m$ is closed in $L_2(\phi_m)$ and the inverse mapping
$$
R_m =(J_m\Phi_mI_m)^{-1}: W_m \to Y_m
$$
satisfies $\|R_m\|\le (1-\ep2^{-m})^{-1} \gamma_p.$

We now introduce the orthogonal projections
$$
P_m : L_2(\phi_m) \to W_m
$$
and consider $Q_m: L_p\to Y_m$ defined to be $Q_m= R_mP_mJ_m\Phi_m$.
For $f\in L_p$ we have
\begin{equation*}
\sum\|Q_mf\|_p^p \le \sum \|R_m\|^p \cdot\|\Phi_mf\|^p_{L_p(\phi_m)}
\le (1-\ep)^{-p}\gamma_p^p\sum \|f1_{A_m}\|_p^p
\le (1-\ep)^{-p}\gamma_p^p\|f\|_p^p,
\end{equation*}
the last inequality following by disjointness of the sets $A_m$.  Since we
already know that $Y=\overline{\sum Y_m}$ is naturally isomorphic to
$(\bigoplus Y_m)_p$, we see that the series $\sum Q_mf$ converges to
an element $Qf$ of $Y$.  Moreover, the operator $Q$ thus defined
satisfies $\|Q\|\le \gamma_p/(1-\ep)$.

To finish, we investigate $\|Q(y)-y\|_p$, when $y= \sum y_k$ with
$y_k\in Y_k$.  If, as before, we write $y_k'=y_k1_{A_k}$ we may note
that $Q_k(y_k)=Q_k(y'_k)$ and $Q_m(y'_k)=0$ for $m\ne k$.  Thus
\begin{align*}
\|Q(y)-y\|_p &= \Big\|\sum_k\Big(\sum_m Q_my_k-y_k\Big)\Big\|_p\\
&= \Big\|\sum_k\sum_{m\ne k} Q_my_k\Big\|_p \qquad\text{[since\
 $Q_ky_k=y_k$]}\\
&= \Big\|\sum_k\sum_m Q_m(y_k-y'_k)\Big\|_p\\
&=\Big\|Q\Big(\sum_ky_k-y'_k\Big)\Big\|_p\\
&\le \|Q\|\sum_k\|y_k-y_k'\|_p
\le \gamma_p(1-\ep)^{-1}\sum 2^{-k}\ep \|y_k\|_p,
\end{align*}
using our estimate for $\|Q\|$ and \eqref{yy'} at the last stage.  We
can now see that for suitable chosen $\ep$, $Q$ may be modified to
give a projection $\tilde Q:L_p\to Y$ with $\|\tilde Q\|\le
(1+\theta)\gamma_p$.
\end{proof}


\section{Quotients and embeddings}\label{S:7}
\subsection{Subspaces of $L_p$ that are quotients of $\ell_p\oplus
\ell_2$}

It was shown in \cite{JO2} that a subspace of $L_p$ ($p>2$) that is
isomorphic to a quotient of a subspace of $\ell_p\oplus \ell_2$ is
in fact isomorphic to a subspace of $\ell_p\oplus \ell_2$.  We can
give an alternative proof of this result by applying the main
theorem of this paper.  Clearly all that is needed is to show that
$\ell_p(\ell_2)$ is not a quotient of a subspace of
$\ell_p\oplus\ell_2$.

We shall prove something more general, namely that $\ell_p(\ell_q)$
is not a quotient of a subspace of $\ell_p\oplus \ell_q$ when
$p,q>1$ and $p\ne q$.  By duality it will be enough to consider the
case $p>q$.  For elements $w=(w_1,w_2)$ of $\ell_p\oplus\ell_q$ we
shall write $\|w\|_p = \|w_1\|_p$, $\|w\|_q = \|w_2\|_q$ and
$\|w\|=\|w\|_p\vee \|w\|_q$.

\begin{lem}
Let $1<q<p<\infty$ and let $W$ be a subspace of $\ell_p\oplus
\ell_q$.  Let $X=\ell_q$, let $Q:W\to X$ be a quotient mapping and
let $\lambda$ be a constant with $0<\lambda<\|Q\|^{-1}$. For every
$M>0$ there is a finite-codimensional subspace $Y$ of $X$ such that,
for $w\in W$ we have
$$
\|w\|\le M,\ Q(w)\in Y,\ \|Q(w)\|=1 \implies \|w\|_q >\lambda.
$$
\end{lem}
\begin{proof}
Suppose otherwise.  We can find a normalized block basis $(x_n)$ in
$X$ and elements $w_n$ of $W$ with $\|w_n\|\le M$, $Q(w_n)=x_n$ and
$\|w_n\|_q \le \lambda$. Taking a subsequence and perturbing
slightly, we may suppose that $w_n = w+w_n'$, where $(w_n')$ is a
block basis in $\ell_p\oplus\ell_q$, satisfying $\|w'_n\|\le M$,
$\|w_n'\|_q\le \lambda$.

Since $Q(w) = \text{w-lim}\,Q(w_n) = 0$, we see that $Q(w_n')=x_n$.
We may now estimate as follows using the fact that the $w_n'$ are
disjointly supported:
$$
\Big\|\sum_{n=1}^N w_n'\Big\| = \Big(\sum_{n=1}^N\|w_n'\|_p^p\Big)^{1/p} \vee
\Big(\sum_{n=1}^N \|w_n'\|_q^q\Big)^{1/q}\le N^{1/p}M \vee N^{1/q}\lambda.
$$
Since the $x_n$ are normalized blocks in $X=\ell_q$ we have
\begin{equation*}
N^{1/q} = \Big\|\sum_{n=1}^N x_n\Big\|
\le \|Q\|\,\Big\|\sum_{n=1}^N w_n'\Big\|
\le M\|Q\|N^{1/p} \vee \lambda \|Q\| N^{1/q}.
\end{equation*}
Since $\lambda\|Q\|<1$, this is impossible once $N$ is large enough.
\end{proof}

\begin{prop}\label{NotQuot}
If $1<q<p<\infty$ then $\ell_p(\ell_q)$ is not a quotient of a
subspace of $\ell_p\oplus \ell_q$.
\end{prop}
\begin{proof}
Suppose, if possible that there exists a quotient operator
$$
\ell_p\oplus\ell_q\supseteq Z\Qarrow X = \Big(\bigoplus_{n\in \nat}
X_n\Big)_p
$$
where $X_n = \ell_q$ for all $n$.  Let $K$ be a constant such that
$Q[KB_Z]\supseteq B_X$, let $\lambda$ be fixed with
$0<\lambda<\|Q\|^{-1}$, choose a natural number $m$ with
$m^{1/q-1/p}>K\lambda^{-1}$, and set $M= 2Km^{1/p}$.

Applying the lemma, we find, for each $n$, a finite-codimensional
subspace $Y_n$ of $X_n$ such that
\begin{equation}\label{atleastlambda}
z\in MB_Z,\ Q(z)\in Y_n,\ \|Q(z)\|=1 \implies \|z\|_q >\lambda.
\end{equation}
For each $n$, let $(e_i^{(n)})$ be a sequence in $Y_n$, 1-equivalent
to the unit vector basis of $\ell_q$. For each $m$-tuple $\mathbf i =
(i_1,i_2,\dots,i_m)\in \nat^m$, let $z(\mathbf i) \in Z$ be chosen
with
$$
Q(z(\mathbf i) = e^{(1)}_{i_1}+e^{(2)}_{i_2}+\cdots+e^{(m)}_{i_m},
$$
and $\|z(\mathbf,,, i)\|\le Km^{1/p}$.

Taking subsequences in each co-ordinate, we may suppose that the
following weak limits exist in $Z$
\begin{align*}
 z(i_1,i_2,\dots,i_{m-1}) &= \text{w-lim}_{i_m\to \infty}\,z(i_1,i_2,\dots,i_{m})\\
&\ \vdots\\
z(i_1,i_2,\dots,i_{j}) &= \text{w-lim}_{i_{j+1}\to \infty}\,z(i_1,i_2,\dots,i_{j+1})\\
&\ \vdots\\
 z(i_1) &= \text{w-lim}_{i_2\to \infty} z(i_1,i_2).
\end{align*}
Notice that, for all $j$ and all $i_1,i_2,\dots,i_j$, the following
hold:
\begin{align*}
Q(z(i_1,\dots,i_j) &= e^{(1)}_{i_1}+\cdots+e^{(j)}_{i_j}\\
\|z(i_1,\dots,i_j) \| &\le Km^{1/p}\\
\|z(i_1,\dots,i_j) -z(i_1,\dots,i_{j-1})\| &\le 2Km^{1/p} = M.
\end{align*}
Since $Q(z(i_1,\dots,i_j) -z(i_1,\dots,i_{j-1}))= e_{i_j}^{(j)}\in
S_{Y_j}$ it must be that
\begin{equation}
\|z(i_1,\dots,i_j) -z(i_1,\dots,i_{j-1})\|_q >\lambda,
\label{lambda}\quad\text{ [by \eqref{atleastlambda}].}
\end{equation}

We shall now choose recursively some special $i_j$ in such a way
that $\|z(i_1,\dots,i_j)\|_q> \lambda j^{1/q}$ for all $j$. Start
with $i_1=1$; since $\|z(i_1)\|\le M$ and $Q(z(i_1))=e^{(1)}_{i_1}$
we certainly have $\|z(i_1)\|_q >\lambda$ by \ref{atleastlambda}.
Since $z(i_1,k)-z(i_1)\to 0$ weakly we can choose $i_2$ such that
$z(i_1,i_2)-z(i_1)$ is essentially disjoint from $z(i_1)$. More
precisely, because of \ref{lambda}, we can ensure that
$$
\|z(i_1,i_2)\|_q=\|z(i_1) + (z(i_1,i_2) - z(i_1))\|_q >
(\lambda^q+\lambda^q)^{1/q} =\lambda 2^{1/q}.
 $$
Continuing in this way, we can indeed choose $i_3,\dots, i_m$ in
such a way that
$$
\|z(i_1,\dots,i_j)\|_q \ge \lambda j^{1/q}.
$$
However, for $j=m$ this yields
$
\lambda m^{1/q} \le Km^{1/p}$,
contradicting our initial choice of $m$.
\end{proof}

\begin{remark} \label{BMdist}
The proof we have just given actually establishes the following
quantitative result: if $Y$ is a quotient of a subspace of
$\ell_p\oplus \ell_q$ then the Banach-Mazur distance $d\big(Y,
\big(\bigoplus_{j=1}^m \ell_q\big)_p\big)$ is at least
$m^{|1/q-1/p|}$.
\end{remark}

\subsection{Uniform bounds for isomorphic embeddings}

As we remarked in the introduction, the Kalton--Werner refinement
\cite{KW} of the result of \cite{JO1} gives an almost isometric
embedding of $X$ into $\ell_p$ when $X$ is a subspace of $L_p$
($p>2$), not containing $\ell_2$.  By contrast, the main result of
the present paper does not have an almost isometric version, and
indeed it is easy to see that there is no constant $K$ (let alone
$K=1+\ep$) such that every subspace of $L_p$ not containing
$\ell_p(\ell_2)$ $K$-embeds in $\ell_p\oplus \ell_2$.  It is enough
to consider spaces $X$ of the form $X=\left(\bigoplus
_{j=1}^m\ell_2\right)_p$. A straightforward argument, or an
application of the more general result mentioned in the remark
above, shows that the Banach--Mazur distance from $X$ to a subspace
of $\ell_p\oplus \ell_2$ is at least $m^{1/2-1/p}$.

If we are looking for a ``uniform'' version of our Main Theorem, it
is perhaps not unreasonable to conjecture the existence of a
constant $K$ such that every subspace of $L_p$ not containing
$\ell_p(\ell_2)$ $K$-embeds in some space of the form $\ell_p\oplus_p
\big(\bigoplus _{j=1}^m\ell_2\big)_p$.  However, no such constant
$M$ exists, as is shown by the following proposition.  The structure
of the space $X$ considered below suggests that if there is some
uniform version of our main result then it will involve independent
sums (see \cite{Als}), rather than, or as well as, $\ell_p$ sums.
The  proof of the next result follows a construction
  due to Alspach and could
be compiled from arguments  in \cite[Chapter 2]{Als}. 
The following is a self contained proof.
\begin{prop}\label{NoUnifEmb}
Let $p>2$. For every $K>0$ there is a subspace $X$ of $L_p$, isomorphic
to $\ell_2$,  such that  for all $m\in\nat$, $X$ is not $K$-isomorphic to a subspace of 
 $\ell_p\oplus_p\bigl(\bigoplus_{l=1}^m \ell_2\bigr)_p$.
\end{prop}
\begin{proof}
Fix a constant $M>1$.  Let $\{v_i, z_{j,k}: i,j,k\in\nat\}$ be a family of independent
random variables in $L_p[0,1]$ with distributions defined as follows:
for $i,j\in\nat$, $z_{i,j}$ is $\mathcal N(0,1)$, while $v_i$ is $\{0,M\}$-valued with
$
\mathbb P[v_i=M] = 1-\mathbb P[v_i=0] = M^{-p/2}.$
We set $x_{i,j}= z_{i,j}\sqrt {v_i}$, noting that
$$
\|x_{i,j}\|_p^p = \mathbb E[v_i^{p/2}|z_{i,j}|^p]=\mathbb
E[v_i^{p/2}]\mathbb E[|z_{i,j}|^p] =\gamma_p^p.
$$
We now define $X_i = [x_{i,j}]_{j\in \mathbb N}$ and $X=
[x_{i,j}]_{i,j\in \mathbb N}$. We start by calculating the norm of a
general element of $X$.

Let $x= \sum_{i,j}c_{i,j}x_{i,j}$.  By independence, and properties
of the normal distribution, the distribution of $x$, {\it
conditional on} $v_1,v_2,v_3,\dots$ is $\mathcal N(0,w)$, where
$w=\sum_{i,j} c_{i,j}^2v_i$.  So
\begin{equation}\label{E:ell2}
\|x\|_p^p = \mathbb E\bigl[\mathbb E[|x|^p\mid v_1,v_2,
\dots]\bigr]
= \gamma_p^p \mathbb E\bigl[ (\sum_i\big(\sum_j c_{i,j}^2\big) v_i)^{p/2}\bigr]=\gamma_p^p\big\|\sum a_i v_i\big\|_{p/2}^{p/2},
\end{equation}
where $a_i=\sum_j c_{i,j}^2$, for $i\in \nat$. Let us first note that \eqref{E:ell2} implies that
  $(x_{i,j})$ is equivalent to the  unit vector basis of $\ell_2$.
Indeed,  Jensen's inequality yields
 $$\big\|\sum a_i v_i\big\|_{p/2}^{p/2}\ge \E^{p/2}\big[\sum a_i v_i]=\big(\sum a_i M^{1-p/2}\big)^{p/2}= 
\big(  M^{1/2-p/4}\big(\sum_{i,j} c_{i,j}^2\big)^{1/2}\big)^p.$$
On the other hands,   letting $\tilde v_i= v_i -\E(v_i)=v_i-M^{1-p/2}$,  the triangle inequality
 in $L_{p/2}$ and the fact that  for some $C<\infty$
 (depending on $M$ and $p$)  the sequence $(\tilde v_i)$, as sequence in $L_{p/2}$, is $C$-equivalent to the unit vector basis in $\ell_2$, implies 
\begin{align*}
 \big\|\sum a_i v_i\big\|_{p/2}&\le M^{1-p/2} \sum a_i  +  \big\|\sum a_i \tilde v_i\big\|_{p/2}\\
                            &\le  M^{1-p/2}\sum a_i+ C \big(\sum a_i^2\big)^{1/2}\le 
  (M^{1-p/2} +C)\sum a_i
\end{align*}
and, thus,  $$\big\|\sum a_i v_i\big\|_{p/2}^{p/2}\le \big((M^{1-p/2} +C)^{1/2}\big(\sum_{i,j} c_{i,j}^2)^{1/2}\big)^p,$$
which finishes the proof of our claim  that  $(x_{i,j})$ is equivalent to the  unit basis of $\ell_2$.

We note two special cases of  \eqref{E:ell2}.  First, if $x=x_i\in X_i$ for some $i$
(thus $c_{i',j}=0$ for all $i'\ne i$ and all $j$), we have
$$
\|x_i\|_p = \gamma_p(\sum_j c_{i,j}^2)^{1/2}.
$$
In particular, $\|x_i\|_p = 1$ if and only if $(\sum_j
c_{i,j}^2)^{1/2} = \gamma_p^{-1}$. Secondly, if $x=n^{-1/2}\sum_{i=1}^n
x_i$, where the $x_i$ are normalized elements of $X_i$,
 \begin{align*}
  \|x\|_p = n^{-1/2}\gamma_p\mathbb E\bigl[ (\sum_{i=1}^n(\sum_j c_{i,j}^2)
v_i)^{p/2}\bigr]^{1/p}
 &= n^{-1/2}\mathbb E\bigl[ (\sum_{i=1}^n v_i)^{p/2}\bigr]^{1/p}
 = \|n^{-1}\sum_{i=1}^n v_i\|_{p/2}^{1/2}.
\end{align*}

 Now, by the weak law of large numbers, $n^{-1}\sum_{i=1}^n v_i$
 converges in probability to the constant $\mathbb
 E[v_1]=M^{1-p/2}$.  Because these averages are uniformly bounded
 (by $M$), the convergence holds also for the $L_{p/2}$-norm.  So as
 $n\to \infty$ we have
 $$
  \Big\|n^{-1}\sum_{i=1}^n v_i\Big\|_{p/2}\to  M^{1-p/2}.
  $$
Summarizing, we can say that if $x_i$ are $L_p$-normalized elements
  of $X_i$ then
  \begin{equation}
\Big \|n^{-1/2}\sum_{i=1}^n x_i\Big\|_p = \Big\|n^{-1}\sum_{i=1}^n v_i\Big\|_{p/2}^{1/2}
  \to M^{(2-p)/4} \label{average-est}
 \text{ 
 as $n\to \infty$.}
\end{equation}

Let $T=(T_\ell)_{\ell=0}^m:X\to
 Y=\ell_p\oplus_p \bigl(\bigoplus_{\ell=1}^m \ell_2\bigr)_p$, with $T_0:X\to\ell_p$  
 and $T_i:X\to\ell_2$, for $\ell=1,2\ldots,m$,
  be an isomorphic embedding. We assume
that $\|T(x)\| \ge \|x\|$ for all $x$ and shall show that $\|T\|\ge
M^{(p-2)/4}$.

We note that, for each $i$, the sequence
$\big(T_0(x_{i,j}))\big)_{j=1}^\infty$ is a weakly null sequence in $\ell_p$.
So by taking vectors of the form
$$
x'_{i,k} = \gamma_p^{-1}k^{-1/2}\sum_{r=1}^k \ x_{i,j_r(k)},
$$
with $j_{k-1}(k-1)<j_1(k)<j_2(k)<\cdots<j_k(k)$, we construct an $L_p$-normalized, weakly
null sequence $(x'_{i,k})_{k=1}^\infty$ in $X_i$ with
$\|T_0(x'_{i,k})\|_p\to 0$ as $k\to \infty$.

 Passing to a subsequence,
we may assume that for all $i\in\nat$ and all $\ell=1,2\ldots m$  the sequence $T_\ell(x'_{i,k}))$ tends to a limit 
$\mu_{i,\ell}$  as $k\to \infty$. Since
$\|T(x'_{i,k})\|\ge 1$ and $\|T_0(x'_{i,k}))\|_p\to 0$, it must be
that $\|\mu_i\|_p\ge 1$, where $\mu_i=(\mu_{i,\ell})_{\ell=1}^m$. Passing to  a subsequence in $i$, we may assume
that $\mu_i$ converges to some $\mu\in \real^m$, as $i\to \infty$, with $\|\mu\|_p\ge 1$.

For $\ell=1,2\ldots m$ and $n\in\nat$ we observe that
\begin{align*}
\lim_{k_1\to \infty} &\lim_{k_{2}\to \infty}\dots\lim_{k_n\to \infty}
 \|n^{-1/2}T_\ell\big(\sum_{i=1}^n x'_{i,k_i}\big)\|_2  \\
&=\lim_{k_1\to \infty} \lim_{k_{2}\to \infty}\dots\lim_{k_{n-1}\to \infty}
n^{-1/2}\Big(\big\|T_\ell\big(\sum_{i=1}^{n-1} x'_{i,k_i}\big)\big\|^2+\mu_{i,\ell}^2 \Big)^{1/2}\\
&=\ldots=n^{-1/2} \Big(\sum_{i=1}^n \mu_{i,\ell}^2 \Big)^{1/2}\equiv \tilde\mu_{n,\ell}.
\end{align*}
Since $\tilde\mu_n\to \mu$, as $n\to \infty$, where $\tilde \mu_n =(\tilde\mu_{n,\ell})_{\ell=1}^m$, we deduce
\begin{equation}
\lim_{n\to \infty}\lim_{k_1\to \infty} \lim_{k_{2}\to
\infty}\dots\lim_{k_n\to \infty} \big\|n^{-1/2}T\big(\sum_{i=1}^n
x'_{i,k_i}\big)\big\|_Y=
\lim_{n\to\infty} \|\tilde\mu_n\|_p
= \|\mu\|_p\ge 1. \label{ge-one}
\end{equation}

On the other hand, as we have already noted above
(\ref{average-est}),
$$
\big\|n^{-1/2}\sum_{i=1}^n x'_{i,k_i}\big\| =   \big\|n^{-1}\sum_{i=1}^n
v_i\big\|_{p/2}^{1/2}\to M^{(2-p)/4},\text{ as $n\to\infty$} ,
$$
Comparing this with  
(\ref{ge-one}), we conclude that $\|T\|\ge M^{(p-2)/4}$ as claimed.
\end{proof}


\section{Concluding Remarks}\label{S:8}        

A natural question remains, namely to characterize when a subspace
$X\subseteq L_p$ $(2<p<\infty)$ embeds into $\ell_p(\ell_2)$. We do
not know the answer. In light of the \cite{JO2} $\ell_p\oplus\ell_2$
quotient result (see paragraph 7.1 above) we ask the following.

\begin{prob}\label{prob8.1}
Let $X\!\subseteq\! L_p$ $(2\!<\!p\!<\!\infty)$. If $X$ is a quotient of
$\ell_p(\ell_2)$ does $X$ embed into $\ell_p(\ell_2)$?
\end{prob}

Extensive study has been made of the $\L_p$ spaces, i.e., the
complemented subspaces of $L_p$ which are not isomorphic to $\ell_2$
(see e.g., \cite{LP} and \cite{LR}). In particular there are
uncountably many such spaces \cite{BRS} and even infinitely many
which embed into $\ell_p$ $(\ell_2)$ \cite{S1}.  Thus it seems that a deeper study 
of the index in \cite{BRS} will be needed  for further progress.
 However some things,
which we now recall, are known.

\begin{thm}\label{thm8.2}
\cite{P} If $Y$ is complemented in $\ell_p$ then $Y$ is isomorphic
to $\ell_p$.
\end{thm}

\begin{thm}\label{thm8.3}
 \cite{JZ}  If $Y$ is a $\L_p$ subspace of $\ell_p$ then $Y$ is isomorphic
to $\ell_p$.
\end{thm}

\begin{thm}\label{thm8.4}
\cite{EW}  If $Y$ is complemented in $\ell_p\oplus\ell_2$ then $Y$
is isomorphic to $\ell_p$, $\ell_2$ or $\ell_p \oplus\ell_2$.
\end{thm}

\begin{thm}\label{thm8.5}
\cite{O} If $Y$ is complemented in $\ell_p(\ell_2)$ then $Y$ is
isomorphic to $\ell_p$, $\ell_2$, $\ell_p\oplus\ell_2$ or $\ell_p$
$(\ell_2)$.
\end{thm}

We recall that $X_p$ is the $\L_p$ discovered by H.~Rosenthal
\cite{R}. For $p>2$, $X_p$ may be defined to be the subspace of
$\ell_p\oplus \ell_2$ spanned by $(e_i + w_i f_i)$, where $(e_i)$
and $(f_i)$ are the unit vector bases of $\ell_p$ and $\ell_2$,
respectively, and  where $w_i\to0$ with $\sum w_i^{2p/p-2} =\infty$.
Since  $\ell_p \oplus\ell_2$ embeds into $X_p$, the subspaces of
$X_p$ and of $\ell_p\oplus \ell_2$ are (up to isomorphism) the same.
For $1<p<2$ the space $X_p$ is defined to be the dual of $X_{p'}$
where $1/p+1/p'=1$.  When restricted to $\mathcal L_p$-spaces, the
results of this paper lead to a dichotomy valid for $1<p<\infty$.

\begin{prop}
Let $Y$ be a $\L_p$-space $(1<p<\infty)$.  Either $Y$ is isomorphic
to a complemented subspace of $X_p$ or $Y$ has a complemented
subspace isomorphic to $\ell_p(\ell_2)$.
\end{prop}
\begin{proof}
For $p>2$ it is shown in \cite{JO2} that a $\L_p$-space which embeds
in $\ell_p\oplus \ell_2$ embeds complementedly in $X_p$.  Combining
this with the main theorem of the present paper gives what we want
for $p>2$.  When $1<p<2$, the space $X_p$ is defined to be the dual
of $X_{p'}$ and so a simple duality argument extends the result to
the full range $1<p<\infty$. \end{proof}

It remains a challenging problem to understand more deeply the
structure of the $\L_p$-subspaces of $X_p$ and $\ell_p\oplus
\ell_2$.

\begin{thm}\label{thm8.6}
 \cite{JO2}  If $Y$ is a $\L_p$ subspace of $\ell_p\oplus\ell_2$ (or $X_p$),
$2<p<\infty$, and $Y$ has an unconditional basis then $Y$ is
isomorphic to $\ell_p$, $\ell_p\oplus\ell_2$ or $X_p$.
\end{thm}

It is known \cite{JRZ} that every $\L_p$ space has a basis but it
remains open if it has an unconditional basis.

\begin{thm}\label{thm8.7}
 \cite{JO2}  If $Y$ is a $\L_p$ subspace of $\ell_p\oplus\ell_2$
$(1<p<2)$ with an unconditional basis then $Y$ is isomorphic to
$\ell_p$ or $\ell_p\oplus\ell_2$.
\end{thm}

So the main open problem for small $\L_p$ spaces is to overcome the
unconditional  basis requirement of \ref{thm8.6} and \ref{thm8.7}.

\begin{prob}\label{prob8.8}
(a) Let $X$ be a $\L_p$ subspace of $\ell_p\oplus\ell_2$
$(2<p<\infty)$. Is $X$ isomorphic to $\ell_p$, $\ell_p \oplus\ell_2$
or $X_p$?

(b) Let $X$ be a $\L_p$ subspace of $\ell_p\oplus \ell_2$
$(1<p<2)$. Is $X$ isomorphic to $\ell_p$ or $\ell_p
\oplus\ell_2$?
\end{prob}

\end{document}